\input ./style/arxiv-general.cfg
\documentclass[aos,MSNbibl,nameyear,rotating,seceqn,dvips]{arximspdf}
\makeatletter
   \@ifpackageloaded{graphicx}{}{\usepackage{graphicx}}
\makeatother
\usepackage{dcolumn}

\doi{10.1214/15-AOS1382}
\volume{44}
\issue{2}
\pubyear{2016}
\firstpage{682}
\lastpage{712}
\docsubty{FLA}

\makeatletter

\def\vfrac#1#2{(#1)/#2}

\def\vafrac#1#2{(#1)/(#2)}
\newcolumntype{d}[1]{D{.}{.}{#1}}
\newcommand{\rrVert}{\Vert}
\newcommand{\llVert}{\Vert}
\newtheorem{lem}{Lemma}

\newproclaim{defn}{Definition}
\newtheorem{prop}{Proposition}
\newproclaim{remark}{Remark}
\newproclaim{exam}{Example}[section]
\newcommand{\bA}{\mathbf{A}}
\newcommand{\bB}{\mathbf{B}}
\newcommand{\bI}{\mathbf{I}}
\newcommand{\bM}{\mathbf{M}}
\newcommand{\bQ}{\mathbf{Q}}
\newcommand{\bx}{\mathbf{x}}
\newcommand{\bsigma}{\bolds{\sigma}}
\newcommand{\bgamma}{\bolds{\gamma}}
\newcommand{\brho}{\bolds{\rho}}
\makeatother

\begin{document}
\begin{frontmatter}

\title{Optimal large-scale quantum state tomography with Pauli
measurements}
\runtitle{Optimal quantum state tomography}

\begin{aug}
\author[A]{\fnms{Tony}~\snm{Cai}\thanksref{m1,T1}\ead[label=e1]{tcai@wharton.upenn.edu}},
\author[B]{\fnms{Donggyu}~\snm{Kim}\thanksref{m2}\ead[label=e2]{kimd@stat.wisc.edu}},
\author[B]{\fnms{Yazhen}~\snm{Wang}\corref{}\thanksref{m2,T2}\ead[label=e3]{yzwang@stat.wisc.edu}},
\author[B]{\fnms{Ming}~\snm{Yuan}\thanksref{m2,T3}\ead[label=e4]{myuan@stat.wisc.edu}}
\and
\author[C]{\fnms{Harrison H.}~\snm{Zhou}\thanksref{m3,T4}\ead[label=e5]{huibin.zhou@yale.edu}}
\runauthor{T. Cai et al.}
\affiliation{University of Pennsylvania\thanksmark{m1}, University of
Wisconsin-Madison\thanksmark{m2}\\ and Yale University\thanksmark{m3}}
\address[A]{T. Cai\\
Department of statistics\\
The Wharton School\\
University of Pennsylvania\\
Philadelphia, Pennsylvania 19104\\
USA\\
\printead{e1}}
\address[B]{D. Kim\\
Y. Wang\\
M. Yuan\\
Department of statistics\\
University of Wisconsin-Madison\\
Madison, Wisconsin 53706\\
USA\\
\printead{e2}\\
\phantom{E-mail:\ }\printead*{e3}\\
\phantom{E-mail:\ }\printead*{e4}}
\address[C]{H. H. Zhou\\
Department of statistics\\
Yale University\\
New Haven, Connecticut 06511\\
USA\\
\printead{e5}}
\end{aug}
\thankstext{T1}{Supported in part by NSF
Grants DMS-12-08982 and DMS-14-03708, and NIH Grant R01 CA127334.}
\thankstext{T2}{Supported in part by NSF Grants DMS-10-5635, DMS-12-65203 and DMS-15-28375.}
\thankstext{T3}{Supported in part by  NSF Career Award DMS-13-21692 and FRG Grant DMS-12-65202.}
\thankstext{T4}{Supported in part  by NSF Grant DMS-12-09191.}

%
\received{\smonth{1} \syear{2015}}
%
\revised{\smonth{8} \syear{2015}}

\begin{abstract}
Quantum state tomography aims to determine the state of a quantum
system as represented by a density matrix.
It is a fundamental task in modern scientific studies involving quantum systems.
In this paper, we study estimation of high-dimensional density matrices
based on 
Pauli measurements. In particular, under appropriate notion of
sparsity, we establish the minimax optimal rates of convergence for
estimation of the density matrix under both the spectral and Frobenius
norm losses; and show how these rates can be achieved by a common
thresholding approach. Numerical performance of the proposed estimator
is also investigated.
\end{abstract}

\end{frontmatter}

\section{Introduction}\label{sec1}

For a range of scientific studies including quantum computation,
quantum information and quantum simulation, an important task is to
learn and engineer quantum systems [\citet{Aspetal05},
\citeauthor{BenCasStr04} (\citeyear{BenCasStr04}, \citeyear{BenCasStr07}),
\citet{Bru12}, \citet{Jon13},\break \citet{Lanetal10}, \citet{NieChu00}, and
Wang (\citeyear{Wan11}, \citeyear{Wan12})]. A~quantum
system is described by its state characterized by a density matrix,
which is a positive semidefinite Hermitian matrix with unit trace.
Determining a quantum state, often referred to as quantum state
tomography, is an important but difficult task [\citet{Alqetal13},
\citet{ArtGilGut05}, \citet{AubButMez09},
\citet{ButGutArt07},
\citet{GutArt07}, \citet{Hafetal05}, \citet{Wan13}, and
\citet{WanXu15}]. It is often inferred by performing measurements on
a large number of identically prepared quantum systems.

More specifically, we describe a quantum spin system by the
$d$-dimensional complex space $\mathbb{C}^d$ and its quantum state by
a complex matrix on $\mathbb{C}^d$. When measuring the quantum system
by performing measurements on some observables which can be represented
by Hermitian matrices, we obtain the measurement outcomes for each
observable, where the measurements take values at random from all
eigenvalues of the observable, with the probability of observing a
particular eigenvalue equal to the trace of the product of the density
matrix and the projection matrix onto the eigenspace corresponding to
the eigenvalue. To handle the up and down states of particles in a
quantum spin system, a common approach is to employ the well-known
Pauli matrices as observables to perform measurements and obtain the
so-called Pauli measurements [\citet{Brietal12}, \citet{Johetal11}, \citet{Liu11},
\citet{SakNap10}, \citet{Sha94}, and Wang (\citeyear{Wan12}, \citeyear{Wan13})].
Since all Pauli matrices have $\pm1$ eigenvalues, Pauli measurements
takes discrete values $1$ and $-1$, and the resulted measurement
distributions can be characterized by binomial distributions. The goal
is to estimate the density matrix based on the Pauli measurements.

Traditional quantum tomography employs classical statistical models and
methods to deduce quantum states from quantum measurements. These
approaches are designed for the setting where the size of a density
matrix is greatly exceeded by the number of quantum measurements, which
is almost never the case even for moderate quantum systems in practice
because the dimension of the density matrix grows exponentially\vspace*{1.5pt} in the
size of the quantum system. For example, the density matrix for 
$b$ spin-$\frac{1}{2}$ quantum systems is of size $2^b\times2^b$. In
this paper, we aim to effectively and efficiently reconstruct the
density matrix for a large-scale quantum system with a relatively
limited number of quantum measurements.

Quantum state tomography is fundamentally connected to the problem of
recovering a high-dimensional matrix based on noisy observations [\citet{Wan13}]. The latter problem arises naturally in many applications in
statistics and machine learning and has attracted considerable recent attention.
When assuming that the unknown matrix of interest is of (approximately)
low-rank, many regularization techniques have been developed. Examples
include \citet{CanRec09},
\citet{CanTao10},
\citeauthor{CanPla09} (\citeyear{CanPla09}, \citeyear{CanPla11}),
\citet{KesMonOh10}, \citet{RecFazPar10},
Bunea, She and Wegkamp (\citeyear{BunSheWeg11}, \citeyear{BunSheWeg12}),
Klopp (\citeyear{Klo11}, \citeyear{Klo12}),
\citet{Kol11}, \citet{KolLouTsy11}, \citet{NegWai11}, \citet{Rec11}, \citet{RohTsy11}, and \citet{CaiZha15}, among many others.
Taking advantage of the low-rank structure of the unknown matrix, these
approaches can often be applied to estimate unknown matrices of high
dimensions. Yet these methods do not fully account for the specific
structure of quantum state tomography. As demonstrated in a pioneering
article, \citet{Groetal10} argued that, when considering quantum
measurements characterized by the Pauli matrices, the density matrix
can often be characterized by the sparsity with respect to the Pauli
basis. Built upon this connection, they suggested a compressed sensing
[\citet{Don06}] strategy for quantum state tomography [\citet{Gro11} and
\citet{Wan13}].
Although promising, their proposal assumes exact measurements, which is
rarely the case in practice, and adopts the constrained nuclear norm
minimization method, which may not be an
appropriate matrix completion approach for estimating a density matrix
with unit trace (or unit nuclear norm). We specifically address such
challenges in the present paper. In particular, we establish the
minimax optimal rates of convergence for the density matrix estimation
under both the spectral and Frobenius norm losses when assuming that
the true density matrix is approximately sparse under the Pauli basis.
Furthermore, we show that these rates could be achieved by carefully
thresholding the coefficients with respect to the Pauli basis. Because
the quantum Pauli measurements are characterized by the binomial
distributions, the convergence rates and minimax lower bounds are
derived by asymptotic analysis with manipulations of binomial
distributions instead of the usual normal distribution based calculations.

The rest of paper proceeds as follows. Section~\ref{sec2} gives some background
on quantum state tomography and introduces a thresholding based density
matrix estimator. Section~\ref{sec3} develops theoretical properties for the
density matrix estimation problem. In particular, the convergence rates
of the proposed density matrix estimator and its minimax optimality
with respect to both the spectral and Frobenius norm losses are
established. Section~\ref{sec4} features a simulation study to
illustrate finite sample performance of the proposed estimators. All
technical proofs are collected in Section~\ref{proofs}.

\section{Quantum state tomography with Pauli measurements}\label{sec2}
In this section, we first review the quantum state and density matrix
and introduce Pauli matrices and Pauli measurements. We also develop
results to describe density matrix representations through Pauli
matrices and characterize the distributions of Pauli measurements via
binomial distribution before introducing a thresholding based density
matrix estimator.

\subsection{Quantum state and measurements}\label{sec21}

For a $d$-dimensional quantum system, we describe its quantum state by
a density matrix $\brho$ on $d$ dimensional complex
space $\mathbb{C}^d$, where density matrix $\brho$ is a $d$ by $d$
complex matrix satisfying (1) Hermitian, that is, $\brho$ is equal to
its conjugate transpose; (2) positive semidefinite; (3)~unit trace,
that is, $\operatorname{tr}(\brho)=1$.

For a quantum system, it is important but difficult to know its quantum
state. Experiments are conducted to perform measurements on the quantum
system and obtain data for studying the quantum system and estimating
its density matrix. In physics literature, quantum state tomography
refers to reconstruction of a quantum state based on measurements for
the quantum systems. Statistically, it is the problem of estimating the
density matrix from the measurements. Common quantum measurements are
on observable $\bM$, which is defined as a Hermitian matrix on
$\mathbb{C}^d$. Assume that the observable $\bM$ has the following
spectral decomposition:
%
\begin{equation}
\label{diagonal} \bM= \sum_{a=1}^r
\lambda_a \bQ_a,
\end{equation}
where $\lambda_a$
are $r$ different real eigenvalues of $\bM$, and $\bQ_a$ are
projections onto the eigenspaces corresponding to $\lambda_a$.
For the quantum system prepared in state $\brho$, we need a
probability space $(\Omega, \mathcal{F}, P)$ to describe
measurement outcomes when performing measurements on the observable
$\bM$.
Denote by $R$ the measurement outcome of $\bM$. According to the
theory of quantum
mechanics, $R$ is a random variable on $(\Omega, \mathcal{F}, P)$ taking
values in $\{\lambda_1, \lambda_2, \ldots, \lambda_r\}$,
with probability distribution given by
%
\begin{equation}
\label{measurement}
P(R = \lambda_a) =\operatorname{tr}(\bQ_a \brho),\qquad a=1,
2, \ldots, r,\qquad E(R)= \operatorname{tr}(\bM\brho).
\end{equation}
We may perform measurements on an observable for a quantum system that
is identically prepared under the state and obtain independent
and identically distributed observations. See \citet{Hol82}, \citet{SakNap10}, and \citet{Wan12}.

\subsection{Pauli measurements and their distributions} \label{SEC-2-2}

The Pauli matrices as observables are widely used in quantum physics
and quantum information science to perform
quantum measurements. 
Let
\begin{eqnarray*}
\bsigma_0 &=& %
\pmatrix{ 1 & 0
\cr
\vspace*{3pt} 0 & 1 },\qquad
\bsigma_1 = %
\pmatrix{ 0 & 1
\cr
\vspace*{3pt} 1 & 0 },\\
\bsigma_2 &=& %
\pmatrix{ 0 & -\sqrt{-1} \vspace*{3pt}
\cr
\sqrt{-1} & 0 }, \qquad \bsigma_3 = %
\pmatrix{ 1 & 0
\cr
\vspace*{3pt} 0 & -1},
\end{eqnarray*}
where $\bsigma_1$, $\bsigma_2$ and $\bsigma_3$ are called the
two-dimensional Pauli matrices.
Tensor products are used to define high-dimensional Pauli matrices. Let
$d=2^b$ for some integer $b$.
We form $b$-fold tensor products of $\bsigma_0$, $\bsigma_1$,
$\bsigma_2$ and $\bsigma_3$ to obtain
$d$ dimensional Pauli matrices
%
\begin{equation}
\label{Pauli-matrix} \bsigma_{\ell_1} \otimes\bsigma_{\ell_2} \otimes\cdots
\otimes \bsigma_{\ell_b},\qquad (\ell_1, \ell_2,\ldots,
\ell_b) \in\{0, 1, 2, 3\}^b.
\end{equation}
%



We identify index $j =1, \ldots, d^2$ with $(\ell_1, \ell_2,\ldots,
\ell_b) \in\{0, 1, 2, 3\}^b$. 
For example, $j=1$ corresponds to $\ell_1= \cdots= \ell_b=0$. With
the index identification we denote
by $\bB_j$ the Pauli matrix $\bsigma_{\ell_1} \otimes\bsigma_{\ell
_2} \otimes\cdots\otimes\bsigma_{\ell_b}$, with $\bB_1=\bI_d$.
We have the following theorem to describe Pauli matrices and represent
a density matrix by Pauli matrices.

\begin{prop} \label{thm-1}
\textup{(i)} Pauli matrices $\bB_2, \ldots, \bB_{d^2}$ are of full
rank and have eigenvalues $\pm1$.
Denote by $\bQ_{j \pm}$ the projections onto the eigen-spaces of $\bB
_j$ corresponding to eigenvalues $\pm1$,
respectively. Then for $j, j^\prime=2, \ldots, d^2$,
\begin{eqnarray*}
&& \operatorname{tr}(\bQ_{j \pm}) = \frac{d}{2}, \qquad\operatorname{tr}(\bB_{j^\prime}
\bQ_{j
\pm}) = \cases{
\displaystyle \pm\frac{d}{2},& \quad$\mbox{if } j=j^\prime,$
\vspace*{3pt}\cr
0,& \quad$\mbox{if }j \neq j^\prime$.}
\end{eqnarray*}
\textup{(ii)}
Denote by $\mathbb{C}^{d\times d}$ the space of all $d$ by $d$ complex
matrices equipped with the Frobenius norm.
All Pauli matrices defined by (\ref{Pauli-matrix}) form an orthogonal
basis for all complex Hermitian matrices.
Given a density matrix $\brho$, we can expand it under the Pauli basis
as follows:
%
\begin{equation}
\label{beta-representation}
\brho= \frac{\bI_d}{d} + \sum_{j=2}^{d^2}
\beta_j \frac{\bB_j}{d},
\end{equation}
where $\beta_j$ are coefficients. 
For $j=2, \ldots, d^2$,
\begin{eqnarray*}
&& \operatorname{tr}(\brho\bQ_{j\pm}) = \frac{1\pm\beta_j}{2}.
\end{eqnarray*}
\end{prop}


Suppose that an experiment is conducted to perform measurements on
Pauli observable $\bB_j$ independently for $n$
quantum systems which are identically prepared in the same quantum
state $\brho$. As $\bB_j$ has eigenvalues $\pm1$,
the Pauli measurements take values $1$ and $-1$, and thus the average
of the $n$ measurements for each $\bB_j$
is a sufficient statistic. Denote by $N_j$ the average of the $n$
measurement outcomes obtained from measuring
$\bB_j$, $j=2,\ldots, d^2$. Our goal is to estimate $\brho$ based on
$N_2, \ldots, N_{d^2}$.

The following proposition provides a simple binomial characterization
for the distributions of $N_j$.

\begin{prop} \label{thm-2}
Suppose that $\brho$ is given by (\ref{beta-representation}). Then
$N_2, \ldots, N_{d^2}$ are independent with
\[
E(N_j)=\beta_j,\qquad\operatorname{Var}(N_j) =
\frac{1-\beta_j^2}{n},
\]
and $n ( N_j + 1)/2$ follows a binomial distribution with $n$ trials
and cell probabilities $\operatorname{tr}(\brho\bQ_{j + })=(1 + \beta_j)/2$,
where $\bQ_{j +}$ denotes the projection onto the eigenspace of $\bB
_j$ corresponding to eigenvalue $1$, and $\beta_j$
is the coefficient of $\bB_j$ in the expansion of $\brho$ in (\ref
{beta-representation}).
\end{prop}

\subsection{Density matrix estimation}\label{sec23}
Since the dimension of a quantum system grows exponentially with its
components such as the number of particles in the system,
the matrix size of $\brho$ tends to be very large even for a moderate
quantum system. We need to impose some structure such as sparsity
on $\brho$ in order to make it consistently estimable. Suppose that
$\brho$ has a sparse representation under the Pauli basis, following
wavelet shrinkage estimation
we construct a density matrix estimator of $\brho$. Assume that
representation (\ref{beta-representation}) is sparse in a sense that
there is only a relatively small number of coefficients $\beta_k$ with
large magnitudes. Formally, we specify sparsity by assuming that
coefficients $\beta_2, \ldots, \beta_{d^2}$ satisfy
%
\begin{equation}
\label{Csparsity}
\sum_{k=2}^{d^2} |
\beta_k|^q \leq\pi_n(d),
\end{equation}
where $0 \leq q <1$, and $\pi_n(d)$ is a deterministic function with
slow growth in $d$
such as $\log d$.

Pauli matrices are used to describe the spins of spin-$\frac{1}{2}$
particles along different directions, and density matrix $\brho$ in
(\ref{beta-representation})
represents a mixture of quantum states with spins along many
directions. Sparsity assumption (\ref{Csparsity}) with $q=0$
indicates the mixed state involving spins along a relatively small
number of directions corresponding to those Pauli matrices with nonzero
$\beta_k$.
The sparsity reduces the complexity of mixed states. Sparse density
matrices often occur in quantum systems where particles have sparse
interactions such as location interactions.
Examples include many quantum systems in quantum information and
quantum computation [\citet{Beretal14}, \citet{Boietal14}, \citet{Brietal12},
\citet{Flaetal12}, \citet{Senetal14}, and Wang (\citeyear{Wan11}, \citeyear{Wan12})].

Since $N_k$ are independent, and $E(N_k)=\beta_k$. We naturally
estimate $\beta_k$ by $N_k$ and
threshold $N_k$ to estimate large $\beta_k$, ignoring small $\beta
_k$, and obtain
%
\begin{eqnarray}
 \hat{\beta}_k & =&  N_k 1\bigl(|N_k|
\geq\varpi\bigr) \quad\mbox{or}
\nonumber
\\[-8pt]
\label{threshold1}
\\[-8pt]
\nonumber
\hat{\beta }_k &=& \operatorname{sign}(N_k)
\bigl(|N_k|- \varpi\bigr)_+, \qquad k = 2, \ldots, d^2,
\end{eqnarray}
and then we use $\hat{\beta}_k$ to construct the following estimator
of $\brho$,
%
\begin{equation}
\label{threshold-estimator1} \hat{\brho} = \frac{\bI_d}{d} + \sum
_{k=2}^{d^2} \hat{\beta}_k
\frac
{\bB_k}{d},
\end{equation}
where the two estimation methods in (\ref{threshold1}) are called hard
and soft thresholding rules, and $\varpi$ is a threshold value which,
we reason below, can be chosen to be
$\varpi= \hbar\sqrt{(4/n) \log d}$ for some constant $\hbar>1$. The
threshold value is designed such that for small $\beta_k$, $N_k$ must be
bounded by threshold $\varpi$ with overwhelming probability, and the
hard and soft thresholding rules select only those $N_k$
with large signal components $\beta_k$.

As $n ( N_k + 1)/2 \sim \operatorname{Bin}(n, (1+\beta_k)/2)$, an application of
Bernstein's inequality leads to
that for any $x>0$,
\begin{eqnarray*}
&& P\bigl( |N_k - \beta_k| \geq x\bigr) 
\leq2
\exp \biggl( - \frac{ n x^2}{ 2( 1 - \beta_k^2 + x /3)} \biggr) \leq2 \exp \biggl( - \frac{n x^2}{2(1 + x /3)}
\biggr),
\end{eqnarray*}
and
\begin{eqnarray*}
&& P \Bigl( \max_{2 \leq k \leq d^2} | N_k -
\beta_k| \leq\varpi \Bigr) 
\\
&&\qquad= \prod
_{k=2}^{d^2} P \bigl( |N_k -
\beta_k | \leq\varpi \bigr)
\\
&&\qquad\geq \biggl[ 1 - 2 \exp \biggl( - \frac{n \varpi^2}{2(1+ \varpi
/3) } \biggr)
\biggr]^{d^2-1} = \bigl[ 1 - 2 d^{-2 \hbar/(1+o(1)) } \bigr]^{d^2-1}
\rightarrow1, 
\end{eqnarray*}
as $d \rightarrow\infty$ and $\varpi\to0$, that is, with
probability tending to one, $|N_k| \leq\varpi$ uniformly for $k=2,
\ldots, d^2$.
Thus, we can select $\varpi= \hbar\sqrt{(4/n) \log d}$ to threshold
$N_k$ and obtain $\hat{\beta}_k$ in (\ref{threshold1}).

\section{Asymptotic theory for the density matrix
estimator}\label{sec3}

\subsection{Convergence rates}\label{sec31}
We fix matrix norm notation for our asymptotic analysis.
Let $\bx=(x_1, \ldots, x_d)^T$ be a $d$-dimensional vector
and $\bA=(A_{ij})$ be a $d$ by $d$ matrix, and
define their $\ell_\alpha$ norms
\[
\| \bx\|_\alpha= \Biggl( \sum_{i=1}^d
|x_i|^\alpha \Biggr)^{1/\alpha}, \qquad\| \bA\|_\alpha=
\sup\bigl\{ \|\bA \bx\|_\alpha, \|\bx\|_\alpha=1 \bigr\},\qquad  1 \leq\alpha\leq
\infty.
\]
Denote by
$\|\bA\|_F = \sqrt{\operatorname{tr}(\bA^\dagger\bA)}$ the Frobenius norm of $\bA$.

For the case of matrix, the $\ell_2$ norm is called the matrix spectral
norm or operator norm. $\| \bA\|_2$ is equal to the square root of the largest
eigenvalue of $\bA \bA^\dagger$,
%
\begin{equation}
\label{ell-1infty-norm} \| \bA\|_1 = \max_{1 \leq j \leq d} \sum
_{i=1}^d |A_{ij}|,\qquad \| \bA
\|_\infty= \max_{1 \leq i \leq d} \sum
_{j=1}^d |A_{ij}|,
\end{equation}
and
%
\begin{equation}
\label{ell-12infty-norm} \|\bA\|_2^2 \leq\| \bA\|_1
\| \bA\|_\infty.
\end{equation}
For a real symmetric or complex Hermitian matrix $\bA$, $\|\bA\|_2$
is equal to the largest
absolute eigenvalue of $\bA$, $\|\bA\|_F$ is the square root of the
sum of squared eigenvalues,
$\|\bA\|_F \leq\sqrt{d}   \|\bA\|_2$, and (\ref
{ell-1infty-norm})--(\ref{ell-12infty-norm}) imply that
$\| \bA\|_2 \leq \| \bA\|_1=\| \bA\|_\infty$.

The following theorem gives the convergence rates for $\hat{\brho}$
under the spectral and Frobenius norms.
%
\begin{thm} \label{thm-3}
Denote by $\Theta$ the class of density matrices satisfying the
sparsity condition (\ref{Csparsity}).
Assume $n^{c_0} \leq d \leq e^{n ^{c_1}} $ for some constants $c_0>0$
and $c_1<1$. For density matrix estimator
$\hat{\brho}$ defined by (\ref{threshold1})--(\ref
{threshold-estimator1}) with threshold $\varpi= \hbar\sqrt{(4/n)
\log d}$
for some constant $\hbar>1$, we have
\begin{eqnarray*}
\sup_{\brho\in\Theta} E\bigl[\| \hat{\brho} - \brho\|_2^2
\bigr] &\leq &  c_2 \pi^2_n(d) \frac{1}{d^2}
\biggl( \frac{ \log d }{n} \biggr)^{1-q},
\\
\sup_{\brho\in\Theta} E\bigl[\| \hat{\brho} - \brho\|_F^2
\bigr] &\leq &  c_3 \pi_n(d) \frac{1}{d} \biggl(
\frac{ \log d }{n} \biggr)^{1-q/2},
\end{eqnarray*}
where $c_2$ and $c_3$ are constants free of $n$ and $d$.
\end{thm}

\begin{remark}\label{rem1}
Theorem~\ref{thm-3} shows that $\hat{\brho}$
achieves the convergence rate $\pi_n(d) d^{-1} (n^{-1} \log
d)^{1-q/2}$ under the squared Frobenius norm loss and the convergence
rate $\pi_n^2(d) d^{-2} (n^{-1} \log d)^{1-q}$ under the squared
spectral norm loss. Both rates will be shown to be optimal in the next section.
Similar to the optimal convergence rates for large covariance and
volatility matrix estimation [\citet{CaiZho12} and \citet{TaoWanZho13}],
the optimal convergence rates here have factors involving $\pi_n(d)$
and $\log d /n$. However, unlike the covariance and volatility matrix
estimation case,
the convergence rates in Theorem~\ref{thm-3} have factors $d^{-1}$ and
$d^{-2}$ for the squared spectral and Frobenius norms, respectively,
and 
go to zero as $d$ approaches to infinity. In particular, the result
implies that MSEs
of the proposed estimator get smaller for large $d$. This is quite
contrary to large covariance and volatility matrix estimation where the
traces are typically
diverge, the optimal convergence rates grow with the logarithm of
matrix size, and the corresponding MSEs increase in matrix size. The new
phenomenon may be due to the unit trace constraint on density matrix
and that the density matrix representation (\ref{beta-representation})
needs a scaling factor $d^{-1}$ to satisfy the constraint.
Also for finite sample $\hat{\brho}$ may not be positive
semidefinite, we may project $\hat{\brho}$ onto the cone formed by
all density matrices under a given matrix norm $\| \cdot\|$, and
obtain a positive semidefinite density matrix estimator $\tilde{\brho
}$. Since the underlying true density matrix $\brho$ is positive
semidefinite with unit trace, and the representation (\ref
{threshold-estimator1}) ensures that $\hat{\brho}$ has unit trace,
the projection implies
$\| \tilde{\brho} - \hat{\brho}\| \leq\| \brho- \hat{\brho}\|$.
Thus, $\|\tilde{\brho} - \brho\| \leq\| \tilde{\brho} - \hat
{\brho}\| + \|\hat{\brho} - \brho\|
\leq2 \| \hat{\brho} - \brho\|$.
Taking $\| \cdot\|$ as the spectral norm or the Frobenius norm and
using Theorem~\ref{thm-3}, we conclude that $\tilde{\brho}$ has the
same convergence rates as $\hat{\brho}$.
\end{remark}

\subsection{Optimality of the density matrix estimator}\label{sec32}

The following theorem establishes a minimax lower bound for estimating
$\brho$ under the spectral norm.
%
\begin{thm} \label{thm-4}
We assume that $\pi_n(d)$ in the sparsity condition (\ref
{Csparsity}) satisfies
%
\begin{equation}
\label{cond-pi} 
\pi_n(d) \leq\aleph
d^{v} ( \log d )^{q/2} /n^{q/2},
\end{equation}
for some constant $\aleph>0$ and $0<v<1/2$.
Then
\[
\inf_{\check{\brho}} \sup_{\brho\in\Theta} E\bigl[\| \check{
\brho} - \brho\|_2^2\bigr] \geq c_4
\pi^2_n(d) \frac
{1}{d^2} \biggl( \frac{ \log d }{n}
\biggr)^{1-q},
\]
where $\check{\brho}$ denotes any estimator of $\brho$ based on
measurement data $N_2, \ldots, N_{d^2}$, and
$c_4$ is a constant free of $n$ and $d$.
\end{thm}

\begin{remark}\label{rem2}
The lower bound in Theorem~\ref{thm-4} matches the
convergence rate of $\hat{\brho}$ under the spectral norm in Theorem~\ref{thm-3},
so we conclude that $\hat{\brho}$ achieves the optimal convergence rate
under the spectral norm.
To establish the minimax lower bound in Theorem~\ref{thm-4}, we
construct a special subclass of density
matrices and then apply Le Cam's lemma. Assumption (\ref{cond-pi}) is
needed to guarantee the positive definiteness of the constructed
matrices as
density matrix candidates and to ensure the boundedness below from zero
for the total variation of related probability distributions in Le
Cam's lemma.
Assumption (\ref{cond-pi}) is reasonable in a sense that if the
right-hand side of (\ref{cond-pi}) is large enough,
(\ref{cond-pi}) will not impose very restrictive condition on $\pi_n(d)$.
We evaluate the dominating factor $n^{-q/2} d^{v}$ on the right-hand
side of (\ref{cond-pi}) for various scenarios.
First, consider $q=0$, the assumption becomes $\pi_n(d) \leq\aleph
  d^{v}$,
$v<1/2$, and so assumption (\ref{cond-pi}) essentially requires
$\pi_n(d)$ grows in $d$ not faster than $d^{1/2}$, which is not
restrictive at all as $\pi_n(d)$ usually grows slowly in $d$.
The asymptotic analysis of high-dimensional statistics usually allows
both $d$ and $n$ go to infinity. Typically,
we may assume $d$ grows polynomially or exponentially in $n$.
If $d$ grows exponentially in $n$, that is, $d \sim\exp(b_0  n)$ for
some $b_0 > 0$, then $n^{q/2}$ is negligible in
comparison with $d^{v}$, and $n^{-q/2} d^{v}$ behavior like $d^{v}$.
The assumption
in this case is not very restrictive. For the case of polynomial
growth, that is, $d \sim n^{b_1}$ for some $b_1>0$, then
$n^{-q/2} d^{v} \sim d^{v - q/(2 b_1)}$.
If $v - q/ (2b_1)>0$, $n^{-q/2} d^{v} $ grows in $d$ like some positive
power of $d$.
Since we may take $v$ arbitrarily close to $1/2$, the positiveness of
$v - q/(2 b_1)$ essentially requires $b_1 > q$, 
which can often be quite realistic given that $q$ is usually very small.
\end{remark}

The theorem below provides a minimax lower bound for estimating $\brho
$ under the Frobenius norm.
%
\begin{thm} \label{thm-5}
We assume that $\pi_n(d)$ in the sparsity condition (\ref
{Csparsity}) satisfies
%
\begin{equation}
\label{cond-pi-0} 
\pi_n(d) \leq
\aleph^\prime d^{v^\prime} /n^{q},
\end{equation}
for some constants $\aleph^\prime>0$ and $0<v^\prime<2$. Then
\[
\inf_{\check{\brho}} \sup_{\brho\in\Theta} E\bigl[\| \check{
\brho} - \brho\|_F^2\bigr] \geq c_5
\pi_n(d) \frac{1}{d} \biggl( \frac{ \log d }{n}
\biggr)^{1-q/2},
\]
where $\check{\brho}$ denotes any estimator of $\brho$ based on
measurement data $N_2, \ldots, N_{d^2}$, and
$c_5$ is a constant free of $n$ and $d$.
\end{thm}

\begin{remark}\label{rem3}
The lower bound in Theorem~\ref{thm-5} matches the
convergence rate of $\hat{\brho}$ under the Frobenius norm in Theorem~\ref{thm-3},
so we conclude that $\hat{\brho}$ achieves the optimal convergence rate
under the Frobenius norm.
Similar to the Remark~\ref{rem2} after Theorem~\ref{thm-4}, we need to apply
Assouad's lemma to establish the minimax lower
bound in Theorem~\ref{thm-5}, and assumption (\ref{cond-pi-0}) is
used to guarantee the positive definiteness of the constructed matrices
as density matrix candidates and to ensure the boundedness below from
zero for the total variation of related probability distributions in
Assouad's lemma. Also the appropriateness of (\ref{cond-pi-0}) is more
relaxed than (\ref{cond-pi}), as $v^\prime<2$ and the right-hand side
of (\ref{cond-pi-0}) has main powers more than the square of that of
(\ref{cond-pi}).
\end{remark}

It is interesting to consider density matrix estimation under a
Schatten norm, where given a matrix $\bA$ of size $d$, we define its
Schatten $s$-norm by
\[
\| \bA\|_{*s} = \Biggl( \sum_{j=1}^d
| \lambda_j |^s \Biggr)^{1/s},
\]
and $\lambda_1, \ldots, \lambda_d$ are the eigenvalues of the square
root of $\bA^\dagger\bA$.
Spectral norm and Frobenius norm are two special cases of the Schatten
$s$-norm with $s=2$ and \mbox{$s=\infty$}, respectively, and the nuclear norm
corresponds to the Schatten $s$-norm with $s=1$. The following result
provides the convergence rate for the proposed thresholding estimator
under the Schatten $s$-norm loss for $1\le s \le\infty$.

\begin{prop}
\label{schatten.prop}
Under the assumptions of Theorem~\ref{thm-3}, the density matrix
estimator $\hat{\brho}$ defined by (\ref{threshold1})--(\ref
{threshold-estimator1}) with threshold $\varpi= \hbar\sqrt{(4/n)
\log d}$ for some constant $\hbar>1$ satisfies
%
\begin{equation}
\label{Schatten-s1}\qquad \sup_{\brho\in\Theta} E\bigl[\| \hat{\brho} - \brho
\|_{*s}^2\bigr] \leq c \bigl[\pi_n(d)\bigr]
^{2-2/\max(s,2)} \frac{1}{d^{2-2/s}} \biggl(\frac
{\log d}{n} \biggr)
^{1-q+q/\max(s,2)}
\end{equation}
for $1\le s \le\infty$, where $c$ is a constant not depending on $n$
and $d$.
\end{prop}

The upper bound in (\ref{Schatten-s1}) matches the minimax convergence
rates for both the spectral norm and Frobenius norm. Moreover, for the
case of the nuclear norm corresponding to the Schatten $s$-norm with
$s=1$, (\ref{Schatten-s1}) leads to an upper bound with the
convergence rate $\pi_n (d)  (\frac{\log d}{n}  )^{1-
q/2}$. We conjecture that the upper bounds in (\ref{Schatten-s1}) are
rate-optimal under the Schatten $s$-norm loss for all $1\le s \le
\infty$. However, establishing a matching lower bound for the general
Schatten norm loss is a difficult task, and we believe that a new
approach is needed for studying minimax density matrix estimation under
the Schatten $s$-norm, particularly the nuclear norm.

\begin{remark}\label{rem4}
The Pauli basis expansion (\ref{beta-representation}) is orthogonal
with respect to the usual Euclidean inner product, and as in the proof
of Lemma~\ref{lem-norm} we have
\[
\| \hat{\brho} - \brho\|_F^2 = \sum
_{k=2}^{d^2} |\hat{\beta}_k -
\beta_k|^2/d,
\]
where $\hat{\beta}$ and $\hat{\brho}$ are threshold estimators of
$\beta$ and $\brho$, respectively. The sparse vector estimation
problem is well studied under the Gaussian or sub-Gaussian noise case
[\citet{DonJoh94} and \citet{Zha12}] and can be used to
recover the minimax result for density matrix estimation under the
Frobenius norm loss, because of orthogonality. In fact, our relatively
simple proof of the minimax results for the Frobenius norm loss is
essentially the same as the sparse vector estimation approach. However,
such an equivalence between sparse density matrix estimation and sparse
vector estimation breaks down for the general Schatten norm loss such
as the commonly used spectral norm and nuclear norm losses.
For the spectral norm loss, Lemma~\ref{lem-norm} in Section~\ref{proofs} provides a sharp upper bound for $E[\| \hat{\brho} - \brho\|
_{2}^2]$ through the $\ell_1$-norm of $(\hat{\beta}_k - \beta_k)$,
and the proof of the minimax lower bound in Theorem~\ref{thm-4} relies
on the property that the spectral norm is determined by the largest
eigenvalue only. Such a special property allows us to reduce the
problem to a simple subproblem and establish the lower bound under the
spectral norm loss. The arguments cannot be applied to the case of the
general Schatten norm loss in particular the nuclear norm loss.
Moreover, instead of directly applying Lemma~\ref{lem-norm} and Remark~\ref{rem5} in Section~\ref{proofs} to derive upper bounds for the general
Schatten norm loss,
we use the obtained sharp upper bounds for the spectral norm and
Frobenius norm losses together with moment inequalities to derive
sharper upper bounds in Proposition~\ref{schatten.prop}. However,
similar lower bounds are not available. Our analysis leads us to
believe that it is not possible to use sparse vector estimation to
recover minimax lower bound results for the general Schatten norm loss
in particular for the spectral norm loss.
\end{remark}

\section{A simulation study}\label{sec4}
A simulation study was conducted to investigate the performance of the
proposed density matrix estimator for the finite
sample. We took $d=32, 64, 128$ and generated a true density matrix
$\brho$ for each case as follows.
$\brho$ has an expansion over the Pauli basis
%
\begin{equation}
\label{rho-simulation} \bolds{\rho} =d^{-1} \Biggl( \mathbf{I}_d
+ \sum_{j=2}^{d^2} \beta_j
\mathbf{B}_{j} \Biggr),
\end{equation}
where $\beta_j = \operatorname{tr}(\brho\bB_j)$, $j=2, \ldots, d^2$.
From $\beta_2, \ldots, \beta_{d^2}$, we randomly selected $[6 \log
d]$ coefficients $\beta_j$ and set the rest of $\beta_j$ to be zero.
We simulated $[6 \log d]$ values independently from a uniform
distribution on $[-0.2,0.2]$, assigned the simulated values at random
to the selected
$\beta_j$, and then constructed $\brho$ from (\ref{rho-simulation}).
The constructed $\brho$ always has unit trace but may not be positive
semi-definite.
The procedure was repeated until we generated a positive semi-definite
$\brho$. We took it as the true density matrix. The simulation
procedure guarantees the obtained $\brho$ is a density matrix and has
a sparse representation under the Pauli basis.


For each true density matrix $\brho$, as described in Section~\ref{SEC-2-2} we simulated data $N_j$ from a binomial distribution with
cell probability $\beta_j$
and the number of cells $n=100, 200, 500, 1000, 2000$.
We constructed coefficient estimators $\hat{\beta}_j$ by (\ref
{threshold1}) and obtained density matrix estimator
$\hat{\brho}$ using (\ref{threshold-estimator1}). The whole
estimation procedure is repeated $200$ times. The density matrix
estimator is measured by
the mean squared errors (MSE), $E\| \hat{\brho}-\brho\|_2^2$ and
$E\| \hat{\brho}-\brho\|_F^2$, that are evaluated by the average
of $\| \hat{\brho}-\brho\|_2^2$ and
$\| \hat{\brho}-\brho\|_F^2$ over $200$ repetitions, respectively.
Three thresholds were used in the simulation study: the universal threshold
$1.01 \sqrt{ 4 \log d /n}$ for all $\beta_j$, the individual
threshold $1.01 \sqrt{ 4 (1-N_j^2) \log d/n}$ for each $\beta_j$, and
the optimal threshold for all
$\beta_j$, which minimizes the computed MSE for each corresponding
hard or soft threshold method. The individual threshold
takes into account the fact in Theorem~\ref{thm-2} that the mean and
variance of $N_j$ are $\beta_j$ and $(1-\beta_j^2)/n$, respectively,
and the variance\vspace*{1pt}
of $N_j$ is estimated by $(1-N_j^2)/n$.
\begin{figure}\vspace*{-6pt}

\includegraphics{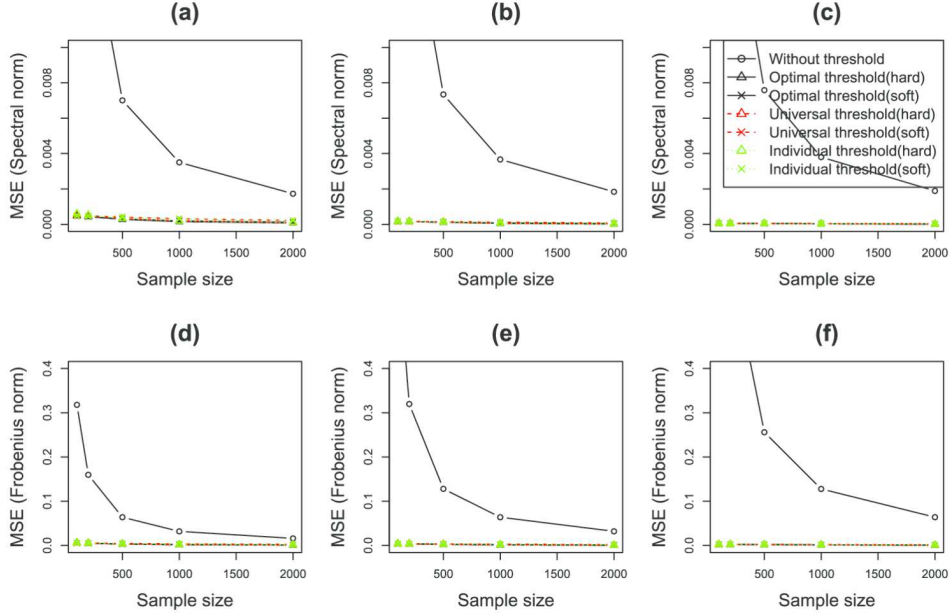}

\caption{The MSE plots against sample size for the proposed density
estimator with hard and soft threshold rules and its corresponding
estimator without thresholding for $d=32, 64, 128$. \textup{(a)}--\textup{(c)} are plots
of MSEs based on the spectral norm for $d=32, 64, 128$, respectively,
and \textup{(d)}--\textup{(f)} are plots of MSEs based on the Frobenius norm for $d=32,
64,128$, respectively.}\vspace*{9pt}
\label{Figure1}
\end{figure}
\begin{figure}

\includegraphics{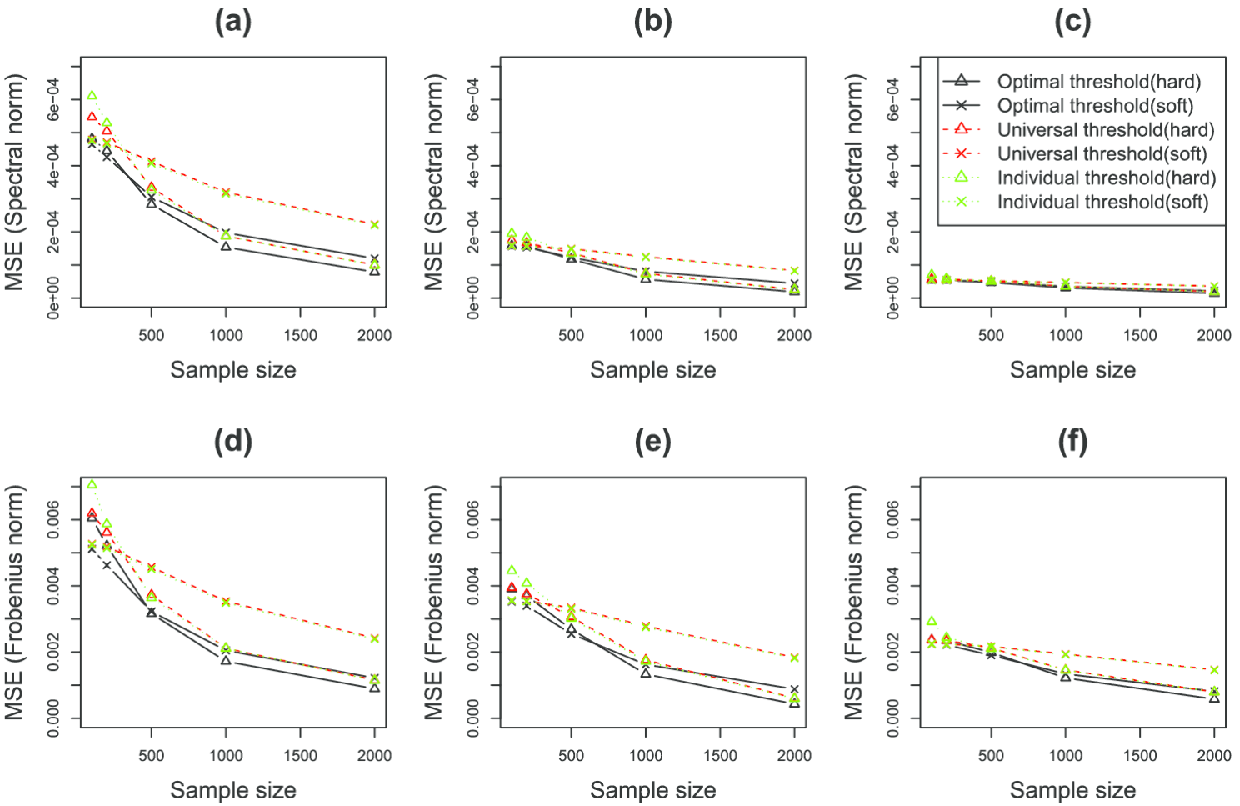}

\caption{The MSE plots against sample size for the proposed density
estimator with hard and soft threshold rules
for $d=32, 64, 128$. \textup{(a)}--\textup{(c)} are plots of MSEs based on the spectral
norm for $d=32, 64, 128$, respectively,
and \textup{(d)}--\textup{(f)} are plots of MSEs based on the Frobenius norm for $d=32,
64, 128$, respectively.}
\label{Figure2}
\end{figure}

\begin{figure}

\includegraphics{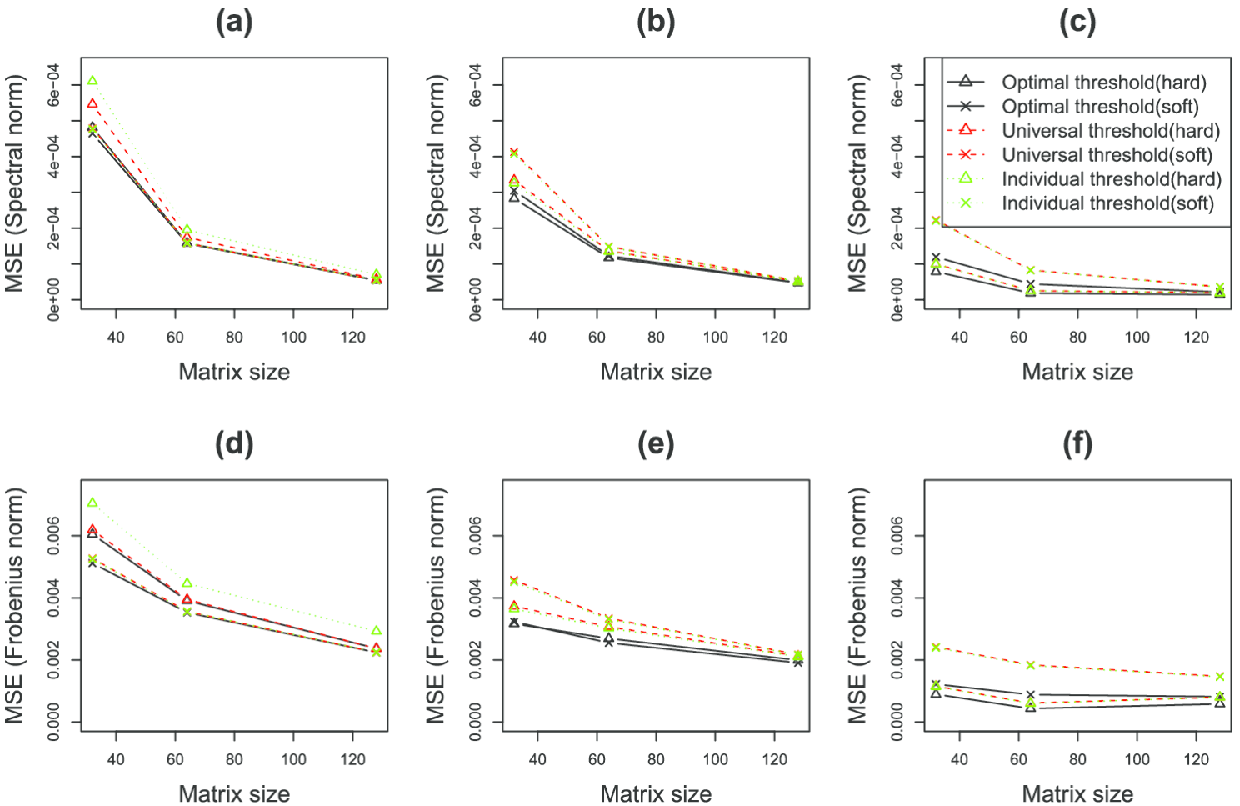}

\caption{The MSE plots against matrix size for the proposed density
estimator with hard and soft threshold rules for $n=100, 500, 2000$.
\textup{(a)}--\textup{(c)} are plots of MSEs based on the spectral norm for $n=100, 500,
2000$, respectively,
and \textup{(d)}--\textup{(f)} are plots of MSEs based on the Frobenius norm for $n=100,
500, 2000$, respectively.}
\label{Figure3}
\end{figure}

\begin{figure}

\includegraphics{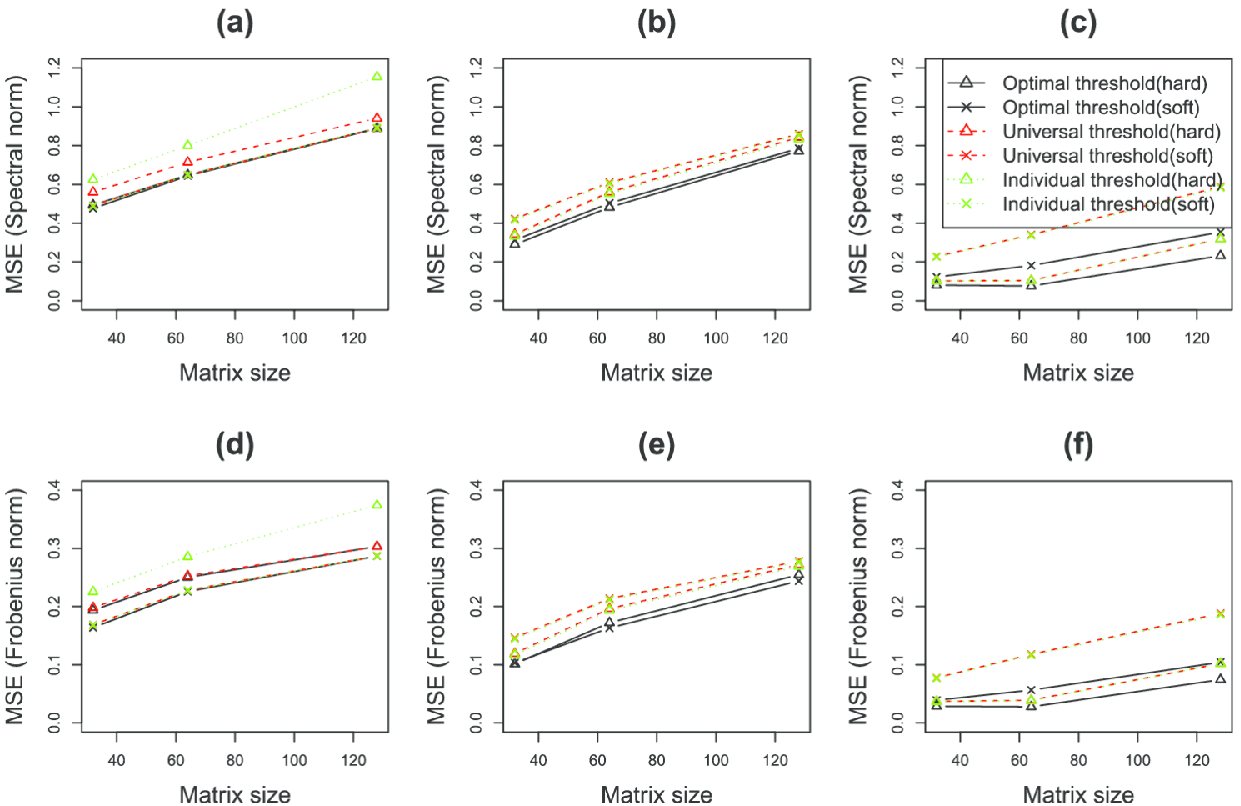}

\caption{The plots of MSEs multiplying by $d$ or $d^2$ against matrix
size $d$ for the proposed density estimator with hard and soft
threshold rules for $n=100, 500, 2000$. \textup{(a)}--\textup{(c)} are plots of $d^2$
times of MSEs based on the spectral norm for $n=100, 500, 2000$, respectively,
and \textup{(d)}--\textup{(f)} are plots of $d$ times of MSEs based on the Frobenius norm
for $n=100, 500, 2000$, respectively.}
\label{Figure4}
\end{figure}

Figures~\ref{Figure1} and \ref{Figure2} plot the MSEs of the density
matrix estimators with hard and soft threshold rules and its
corresponding density matrix estimator without thresholding [i.e.,
$\beta_j$ are estimated by $N_j$ in (\ref{threshold-estimator1})]
against the sample size $n$ for different matrix size $d$, and Figures~\ref{Figure3} and \ref{Figure4} plot their MSEs against matrix size
$d$ for different sample size. The numerical values of the MSEs are
reported in Table~\ref{Table1}.
Figures~\ref{Figure1} and \ref{Figure2} show that the MSEs usually
decrease in sample size $n$, and the thresholding density matrix
estimators enjoy
superior performances than that the density matrix estimator without
thresholding even for $n=2000$; while all threshold rules and threshold
values yield thresholding\vadjust{\goodbreak} density matrix estimators with very close
MSEs, the soft threshold rule with individual and universal threshold
values produce larger MSEs than others for larger sample size such as
$n=1000, 2000$ and the soft threshold rule tends to give somewhat
better performance than the hard threshold rule for smaller
sample size like $n=100, 200$. Figures~\ref{Figure3} and \ref
{Figure4} demonstrate that while the MSEs of all thresholding density
matrix estimators decrease in the matrix size $d$, but if we rescale
the MSEs by multiplying it with $d^2$ for the spectral norm case and
$d$ for the Frobenius norm case, the rescaled MSEs slowly increase in
matrix size $d$. The simulation results largely confirm the theoretical
findings discussed in Remark~\ref{rem1}.

\section{Proofs}\label{proofs}
Let $p=d^2$. Denote by $C$'s generic constants whose values are free of
$n$ and $p$ and
may change from appearance to appearance. Let $u \vee v$ and $u \wedge
v$ be the
maximum and minimum of $u$ and $v$, respectively. For two sequences $u_{n,p}$
and $v_{n,p}$, we write $u_{n,p} \sim v_{n,p}$ if $u_{n,p}/v_{n,p}
\rightarrow1$ as $n, p \rightarrow\infty$,
and write $u_{n,p} \asymp v_{n,p}$ if there exist positive constants
$C_1$ and $C_2$ free of $n$ and $p$ such that
$C_1 \leq u_{n,p}/v_{n,p} \leq C_2$. Without confusion
we may write $\pi_n(d)$ as $\pi_n(p)$.

\subsection{Proofs of Propositions \texorpdfstring{\protect\ref{thm-1}}{1} and
\texorpdfstring{\protect\ref{thm-2}}{2}}\label{sec51}
\mbox{}
\begin{pf*}{Proof of Proposition \protect\ref{thm-1}}
In two dimensions, Pauli matrices satisfy $\operatorname{tr}(\bsigma_0)=2$, and
$\operatorname{tr}(\bsigma_1)=\operatorname{tr}(\bsigma_2)=\operatorname{tr}(\bsigma_3)=0$,
$\bsigma_1, \bsigma_2, \bsigma_3$ have eigenvalues $\pm1$, the
square of a Pauli matrix is equal to the identity matrix, and
the multiplications of
any two Pauli matrices are equal to the third Pauli matrix multiplying
by $\sqrt{-1}$, for example,
$\bsigma_1 \bsigma_2 = \sqrt{-1}\bsigma_3$, $\bsigma_2 \bsigma_3
= \sqrt{-1}\bsigma_1$, and $\bsigma_3 \bsigma_1 = \sqrt{-1}\bsigma_2$.

For $j=2, \ldots, p$, consider
$\bB_j = \bsigma_{\ell_1} \otimes\bsigma_{\ell_2} \otimes\cdots
\otimes\bsigma_{\ell_b}$.
$\operatorname{tr}(\bB_{j})=\operatorname{tr}(\bsigma_{\ell_1}) 
\times\break \operatorname{tr}(\bsigma_{\ell_2})  \cdots
\operatorname{tr}(\bsigma_{\ell_b})=0$, and
$\bB_j$ has eigenvalues $\pm1$, $\bB_j^2 = \bI_d$.

For $j,j^\prime=2, \ldots, p$, $j\neq j^\prime$,
$\bB_j= \bsigma_{\ell_1} \otimes\bsigma_{\ell_2} \otimes\cdots
\otimes\bsigma_{\ell_b}$ and
$\bB_{j^\prime}= \bsigma_{\ell^\prime_1} \otimes\bsigma_{\ell
^\prime_2} \otimes\cdots\otimes\bsigma_{\ell^\prime_b}$,
\[
\bB_j \bB_{j^\prime} = [\bsigma_{\ell_1}
\bsigma_{\ell^\prime
_1}]\otimes[\bsigma_{\ell_2}\bsigma_{\ell^\prime_2}]
\otimes \cdots\otimes[\bsigma_{\ell_b} \bsigma_{\ell^\prime_b}],
\]
is equal to a $d$ dimensional Pauli matrix multiplying by $(\sqrt
{-1})^b$, which has zero trace. Thus,
$\operatorname{tr}(\bB_j \bB_{j^\prime}) = 0$, that is, $\bB_j$ and $\bB
_{j^\prime}$ are orthogonal,
and $\bB_1, \ldots, \bB_p$ form an orthogonal basis.
$\operatorname{tr}(\brho\bB_j/d)= \beta_k \operatorname{tr}(\bB_j^2)/d= \beta_k$. In particular
$\bB_1=\bI_d$, and $\beta_1= \operatorname{tr}(\brho\bB_1) = \operatorname{tr}(\brho)=1$.
\begin{sidewaystable}
\tablewidth=\textwidth
\caption{MSEs based on spectral and Frobenius norms of the density
estimator defined by (\protect\ref{threshold1}) and (\protect\ref
{threshold-estimator1}) and its corresponding density matrix estimator
without thresholding, and threshold values used for $d=32, 64, 128$,
and $n=100, 200, 500, 1000, 2000$}
\label{Table1}
\begin{tabular*}{\textwidth}{@{\extracolsep{4in minus 4in}}lcd{4.3}d{1.3}d{1.3}d{1.3}d{1.3}d{1.3}d{1.3}d{2.3}d{2.3}d{2.3}@{}}
\hline
&
& \multicolumn{7}{c} {\textbf{MSE} \textbf{(Spectral norm)} $\bolds{{\times}10^{4}}$} & \\[-4pt]
&& \multicolumn{7}{l}{\hrulefill}& \\
&&  \multicolumn{1}{c}{\textbf{Without} } &
\multicolumn{2}{c}{\textbf{Optimal}} &
\multicolumn{2}{c}{\textbf{Universal}} &
\multicolumn{2}{c}{\textbf{Individual}}
\\[-18pt]
&&&&&&&&&\multicolumn{3}{c@{}}{\textbf{Threshold value} $\bolds{(\varpi)}$
$\bolds{{\times}10^{2}}$}      \\[-4pt]
&&&&&&&&&\multicolumn{3}{l@{}}{\hrulefill}  \\
&&  \multicolumn{1}{c}{\textbf{threshold}} &
\multicolumn{2}{c}{\textbf{threshold}} &
\multicolumn{2}{c}{\textbf{threshold}} &
\multicolumn{2}{c}{\textbf{threshold}}
& \multicolumn{1}{c}{\textbf{Universal}} & \multicolumn{2}{c@{}}{\textbf{Optimal}}
\\[-4pt]
&& \multicolumn{1}{l}{\hrulefill} & \multicolumn{2}{l}{\hrulefill} &\multicolumn{2}{l}{\hrulefill}& \multicolumn{2}{l}{\hrulefill}& \multicolumn{1}{l}{\hrulefill}&
\multicolumn{2}{l@{}}{\hrulefill}\\
\multicolumn{1}{@{}l}{$\bolds{d}$} &
\multicolumn{1}{c}{$\bolds{n}$} &
\multicolumn{1}{c}{\textbf{Density estimator}} &
\multicolumn{1}{c}{\textbf{Hard}} &
\multicolumn{1}{c}{\textbf{Soft}} &
\multicolumn{1}{c}{\textbf{Hard}} &
\multicolumn{1}{c}{\textbf{Soft}} &
\multicolumn{1}{c}{\textbf{Hard}} &
\multicolumn{1}{c}{\textbf{Soft}} &
\multicolumn{1}{c}{\textbf{Universal}} &
\multicolumn{1}{c}{\textbf{Hard}} &
\multicolumn{1}{c@{}}{\textbf{Soft}} \\
\hline
\phantom{0}32 & \phantom{0}100 & 348.544 & 4.816 & 4.648 & 5.468 & 4.790 & 6.104 & 4.762
& 24.782 & 15.180 & 0.619 \\
& \phantom{0}200 & 175.034 & 4.449 & 4.257 & 5.043 & 4.708 & 5.293 & 4.667 &
17.524 & 7.739 & 0.562 \\
& \phantom{0}500 & 70.069 & 2.831 & 3.054 & 3.344 & 4.130 & 3.260 & 4.071 &
11.083 & 2.397 & 0.373 \\
& 1000 & 35.028 & 1.537 & 1.974 & 1.875 & 3.201 & 1.875 & 3.155 &
7.837 & 1.099 & 0.212 \\
& 2000 & 17.307 & 0.785 & 1.195 & 1.001 & 2.230 & 0.989 & 2.200 &
5.541 & 0.551 & 0.116 \\[3pt]
\phantom{0}64 & \phantom{0}100 & 368.842 & 1.583 & 1.572 & 1.744 & 1.583 & 1.954 & 1.586
& 27.148 & 16.660 & 0.395 \\
& \phantom{0}200 & 183.050 & 1.565 & 1.534 & 1.669 & 1.575 & 1.833 & 1.571 &
19.196 & 9.252 & 0.376 \\
& \phantom{0}500 & 73.399 & 1.175 & 1.228 & 1.367 & 1.490 & 1.347 & 1.476 &
12.141 & 2.900 & 0.307 \\
& 1000 & 36.692 & 0.566 & 0.807 & 0.747 & 1.249 & 0.722 & 1.233 &
8.585 & 1.308 & 0.177 \\
& 2000 & 18.402 & 0.186 & 0.443 & 0.255 & 0.832 & 0.251 & 0.820 &
6.070 & 0.657 & 0.061 \\[3pt]
128 & \phantom{0}100 & 381.032 & 0.543 & 0.542 & 0.574 & 0.543 & 0.705 & 0.545
& 29.323 & 17.500 & 0.237 \\
& \phantom{0}200 & 190.113 & 0.541 & 0.539 & 0.570 & 0.542 & 0.594 & 0.542 &
20.734 & 10.246 & 0.235 \\
& \phantom{0}500 & 75.824 & 0.471 & 0.480 & 0.514 & 0.525 & 0.509 & 0.522 &
13.114 & 3.547 & 0.213 \\
& 1000 & 38.010 & 0.309 & 0.350 & 0.355 & 0.470 & 0.354 & 0.466 &
9.273 & 1.613 & 0.146 \\
& 2000 & 18.907 & 0.142 & 0.216 & 0.194 & 0.359 & 0.194 & 0.356 &
6.557 & 0.725 & 0.080 \\
\hline
\end{tabular*}
\end{sidewaystable}
\setcounter{table}{0}
\begin{sidewaystable}
\tablewidth=\textwidth
\caption{(Continued)}
\begin{tabular*}{\textwidth}{@{\extracolsep{4in minus 4in}}lcd{4.3}d{1.3}d{1.3}d{1.3}d{1.3}d{1.3}d{1.3}d{2.3}d{2.3}d{2.3}@{}}
\hline
&
& \multicolumn{7}{c} {\textbf{MSE} \textbf{(Frobenius norm)} $\bolds{{\times}10^{3}}$} & \\[-4pt]
&& \multicolumn{7}{l}{\hrulefill}& \\
&&  \multicolumn{1}{c}{\textbf{Without} } &
\multicolumn{2}{c}{\textbf{Optimal}} &
\multicolumn{2}{c}{\textbf{Universal}} &
\multicolumn{2}{c}{\textbf{Individual}}
\\[-18pt]
&&&&&&&&&\multicolumn{3}{c@{}}{\textbf{Threshold value} $\bolds{(\varpi)}$
$\bolds{{\times}10^{2}}$}      \\[-4pt]
&&&&&&&&&\multicolumn{3}{l@{}}{\hrulefill}  \\
&&  \multicolumn{1}{c}{\textbf{threshold}} &
\multicolumn{2}{c}{\textbf{threshold}} &
\multicolumn{2}{c}{\textbf{threshold}} &
\multicolumn{2}{c}{\textbf{threshold}}
& \multicolumn{1}{c}{\textbf{Universal}} & \multicolumn{2}{c@{}}{\textbf{Optimal}}
\\[-4pt]
&& \multicolumn{1}{l}{\hrulefill} & \multicolumn{2}{l}{\hrulefill} &\multicolumn{2}{l}{\hrulefill}& \multicolumn{2}{l}{\hrulefill}& \multicolumn{1}{l}{\hrulefill}&
\multicolumn{2}{l@{}}{\hrulefill}\\
\multicolumn{1}{@{}l}{$\bolds{d}$} &
\multicolumn{1}{c}{$\bolds{n}$} &
\multicolumn{1}{c}{\textbf{Density estimator}} &
\multicolumn{1}{c}{\textbf{Hard}} &
\multicolumn{1}{c}{\textbf{Soft}} &
\multicolumn{1}{c}{\textbf{Hard}} &
\multicolumn{1}{c}{\textbf{Soft}} &
\multicolumn{1}{c}{\textbf{Hard}} &
\multicolumn{1}{c}{\textbf{Soft}} &
\multicolumn{1}{c}{\textbf{Universal}} &
\multicolumn{1}{c}{\textbf{Hard}} &
\multicolumn{1}{c@{}}{\textbf{Soft}} \\
\hline
\phantom{0}32 & \phantom{0}100 & 317.873 & 6.052 & 5.119 & 6.195 & 5.274 & 7.050 & 5.246
& 24.782 & 11.004 & 9.936 \\
& \phantom{0}200 & 159.679 & 5.217 & 4.629 & 5.616 & 5.187 & 5.874 & 5.143 &
17.524 & 5.681 & 3.771 \\
& \phantom{0}500 & 63.823 & 3.165 & 3.229 & 3.732 & 4.575 & 3.642 & 4.512 &
11.083 & 2.286 & 0.954 \\
& 1000 & 31.856 & 1.722 & 2.053 & 2.119 & 3.540 & 2.119 & 3.492 &
7.837 & 1.100 & 0.401 \\
& 2000 & 15.967 & 0.894 & 1.219 & 1.155 & 2.424 & 1.141 & 2.394 &
5.541 & 0.546 & 0.182 \\[3pt]
\phantom{0}64 & \phantom{0}100 & 641.437 & 3.909 & 3.528 & 3.951 & 3.563 & 4.463 & 3.562
& 27.148 & 13.719 & 13.234 \\
& \phantom{0}200 & 319.720 & 3.706 & 3.401 & 3.755 & 3.548 & 4.082 & 3.536 &
19.196 & 7.042 & 5.515 \\
& \phantom{0}500 & 127.958 & 2.691 & 2.551 & 3.069 & 3.342 & 3.023 & 3.309 &
12.141 & 2.800 & 1.275 \\
& 1000 & 63.845 & 1.335 & 1.628 & 1.765 & 2.791 & 1.717 & 2.756 &
8.585 & 1.277 & 0.548 \\
& 2000 & 31.952 & 0.433 & 0.882 & 0.610 & 1.842 & 0.596 & 1.817 &
6.070 & 0.647 & 0.258 \\[3pt]
128 & \phantom{0}100 & 1283.182 & 2.370 & 2.240 & 2.370 & 2.242 & 2.924 & 2.245
& 29.323 & 15.989 & 16.128 \\
& \phantom{0}200 & 639.556 & 2.349 & 2.219 & 2.354 & 2.238 & 2.444 & 2.238 &
20.734 & 8.218 & 7.799 \\
& \phantom{0}500 & 255.954 & 1.990 & 1.906 & 2.125 & 2.172 & 2.102 & 2.160 &
13.114 & 3.355 & 1.773 \\
& 1000 & 127.714 & 1.221 & 1.341 & 1.463 & 1.943 & 1.448 & 1.924 &
9.273 & 1.546 & 0.729 \\
& 2000 & 63.921 & 0.581 & 0.815 & 0.798 & 1.471 & 0.798 & 1.456 &
6.557 & 0.719 & 0.327 \\
\hline
\end{tabular*}
\end{sidewaystable}

Denote by $\bQ_{j \pm}$ the projections onto the eigenspaces
corresponding to eigenvalues $\pm1$, respectively.
Then for $j =2, \ldots, p$,
\begin{eqnarray*}
\bB_{j} &=& \bQ_{j +} - \bQ_{j -},\qquad
\bB_j^2=\bQ_{j +} + \bQ _{j -}=
\bI_d, \qquad\bB_{j} \bQ_{j \pm} = \pm
\bQ_{j \pm}^2= \pm \bQ_{j \pm},
\nonumber
\\
0 &=& \operatorname{tr}(\bB_{j}) = \textup{tr}(\bQ_{j +}) - \textup{tr}(
\bQ_{j -}), \qquad d = \operatorname{tr}(\bI _d) = \operatorname{tr}(\bQ_{j +}) + \operatorname{tr}(
\bQ_{j -}),
\end{eqnarray*}
and solving the equations we get
%
\begin{equation}
\label{eq1} \operatorname{tr}(\bQ_{j \pm}) = d/2,\qquad \operatorname{tr}(\bB_{j}
\bQ_{j \pm}) = \pm \operatorname{tr}(\bQ _{j \pm}) = \pm d/2.
\end{equation}
For $j \neq j^\prime$, $j, j^\prime=2, \ldots, p$, $\bB_{j}$ and
$\bB_{j^\prime}$ are orthogonal,
\[
0=\operatorname{tr}(\bB_{j^\prime} \bB_{j}) = \operatorname{tr}(\bB_{j^\prime}
\bQ_{j+}) - \operatorname{tr}(\bB_{j^\prime} \bQ_{j -}),
\]
and
\begin{eqnarray*}
\bB_{j^\prime} \bQ_{j+} + \bB_{j^\prime}
\bQ_{j -} &=& \bB _{j^\prime} (\bQ_{j+} +
\bQ_{j -}) = \bB_{j^\prime},
\\
\operatorname{tr}(\bB_{j^\prime}
\bQ_{j+}) + \operatorname{tr}(\bB_{j^\prime} \bQ_{j -}) &=& \operatorname{tr}(
\bB_{j^\prime})=0,
\end{eqnarray*}
which imply
\begin{eqnarray*}
\label{1}
&& \operatorname{tr}(\bB_{j^\prime} \bQ_{j \pm}) = 0,\qquad  j \neq
j^\prime, j, j^\prime =2, \ldots, p.
\end{eqnarray*}

For a density matrix $\brho$ with representation (\ref
{beta-representation}) under the Pauli basis (\ref{Pauli-matrix}),
from~(\ref{eq1}) we have $\operatorname{tr}(\bQ_{k \pm})=d/2$ and $\operatorname{tr}(\bB_{k} \bQ
_{k \pm}) = \pm d/2$, and thus
%
\begin{eqnarray}
\operatorname{tr}( \brho\bQ_{k \pm}) &=& \frac{1}{d} \operatorname{tr}(
\bQ_{k \pm}) + \sum_{j=2}^p
\frac{\beta_j}{d} \operatorname{tr}(\bB_{j} \bQ_{k \pm} )
\nonumber
\\[-8pt]
\label{tr-Qrho2}
\\[-8pt]
\nonumber
&= &\frac{1}{2}
+ \frac{\beta_k}{d} \operatorname{tr}(\bB_{k} \bQ_{k \pm} ) =
\frac{1 \pm\beta_k}{2}.
\end{eqnarray}
\upqed\end{pf*}

\begin{pf*}{Proof of Proposition~\protect\ref{thm-2}}
We perform measurements on each Pauli observable $\bB_k$ independently
for $n$ quantum systems that are identically prepared under state
$\brho$.
Denote by $R_{k1}, \ldots, R_{kn}$ the $n$ measurement outcomes for
measuring $\bB_k$, $k=2,\ldots, p$.
\begin{equation}
\label{tomography3} N_k=(R_{k1} + \cdots+
R_{kn})/n,
\end{equation}
$R_{k\ell}$, $k=2, \ldots, p$, $\ell=1, \ldots, n$, are
independent, and take values $\pm1$, with distributions given by
%
\begin{eqnarray}
\label{tomography4}\label{tomography5}
&& P(R_{k\ell}= \pm1) = \operatorname{tr}(\brho\bQ_{k \pm} ),\qquad k=2,
\ldots,p, \ell=1, \ldots, n.
\end{eqnarray}

As random variables $R_{k1}, \ldots, R_{kn}$ are i.i.d. and take
eigenvalues $\pm1$,
$n ( N_k + 1)/2= \sum_{\ell=1}^n (R_{k\ell} + 1)/2$ is equal to the
total number of random variables $R_{k1}, \ldots, R_{kn}$
taking eigenvalue $1$, and thus $n ( N_k + 1)/2$ follows a binomial
distribution with $n$ trials and cell probability
$P(R_{k1}= 1) = \operatorname{tr}(\brho\bQ_{k +} )$. From (\ref{tomography3})--(\ref
{tomography4}) and Proposition~\ref{thm-1}, we have
for $k=2, \ldots, p$,
\begin{eqnarray*}
\operatorname{tr}(\brho\bQ_{k +} ) &=& \frac{1+\beta_k}{2},\qquad E(N_k) =
E(R_{k1}) = \operatorname{tr}(\brho\bB_k) = \beta_k \operatorname{tr}
\bigl(\bB_k^2\bigr)/d = \beta _k,
\\
\operatorname{Var}(N_k) &=& \frac{1-\beta_k^2}{n}.
\end{eqnarray*}
\upqed\end{pf*}

%

\subsection{Proof of Theorem \texorpdfstring{\protect\ref{thm-3}}{1}: Upper bound} \label{upper}

\begin{lem} \label{lem-sparse}
If $\beta_j$ satisfy sparsity condition (\ref{Csparsity}), then for
any $a$,
\begin{eqnarray*}
\sum_{j=2}^p |\beta_j| 1\bigl(|
\beta_j| \leq a \varpi\bigr) &\leq &  a^{1-q} \pi _n(p)
\varpi^{1-q},
\\
\sum_{j=2}^p 1\bigl(|\beta_j|
\geq a \varpi\bigr) &\leq &  a^{-q} \pi_n(p) \varpi ^{-q}.
\end{eqnarray*}
\end{lem}
\begin{pf}
Simple algebraic manipulation shows
\begin{eqnarray*}
\sum_{j=2}^p |\beta_j|
1\bigl(|\beta_j| \leq a \varpi\bigr)  &\leq & (a \varpi )^{1-q} \sum
_{j=2}^p |\beta_j|^q
1\bigl(|\beta_j| \leq a \varpi\bigr)
\\
&\leq &  a^{1-q} \pi_n(p) \varpi^{1-q},
\end{eqnarray*}
and
\begin{eqnarray*}
 \sum_{j=2}^p 1\bigl(|\beta_j|
\geq a \varpi\bigr) &\leq & \sum_{j=2}^p \bigl[|
\beta_j|/ (a \varpi)\bigr]^q 1\bigl(|\beta_j|
\geq a \varpi\bigr)
\\
&\leq & (a \varpi)^{-q} \sum_{j=2}^p
|\beta_j|^q \leq a^{-q} \pi _n(p)
\varpi^{-q}.
\end{eqnarray*}
\upqed\end{pf}

\begin{lem} \label{lem-tail}
With $\varpi= \hbar n^{-1/2} \sqrt{2 \log p}$ for some positive
constant $\hbar$, we have for any $a \neq1$,
\begin{eqnarray*}
P\bigl( N_j - \beta_j \leq- |a-1| \varpi\bigr) &\leq & 2
p^{-\vafrac{\hbar^2
|a-1|^2}{1+o(1)}},
\\
P\bigl( N_j - \beta_j \geq|a-1| \varpi\bigr)
& \leq & 2 p^{-\vafrac{\hbar^2 |a-1|^2}{1+o(1)}}.
\end{eqnarray*}
\end{lem}
\begin{pf}
From Proposition~\ref{thm-2} and (\ref{tomography3})--(\ref{tomography4}),
we have that $N_j$ is the average of
$R_{j1}, \ldots, R_{jn}$, which are i.i.d. random variables taking
values $\pm1$,
$P( R_{j1}= \pm1) = (1 \pm\beta_j)/2$, $E(R_{j1})=\beta_j$ and
$\operatorname{Var}(R_{j 1})= 1 - \beta_j^2$.
Applying Bernstein's inequality, we obtain for any $x>0$,
\begin{eqnarray*}
&& P\bigl( |N_j - \beta_j| \geq x\bigr) \leq2 \exp \biggl( -
\frac{ n x^2}{ 2(
1 - \beta_j^2 + x /3)} \biggr) \leq 2 \exp \biggl( - \frac{n x^2}{2(1 + x /3)} \biggr).
\end{eqnarray*}
Both $P( N_j - \beta_j \leq- |a-1| \varpi)$ and $P( N_j - \beta_j
\geq|a-1| \varpi)$ are less than
$P( |N_j - \beta_j| \geq|a-1| \varpi)$, which is bounded by
\begin{eqnarray*}
 2 \exp \biggl( - \frac{ n |a-1|^2 \varpi^2}{ 2(1+ |a-1| \varpi/3)
} \biggr) &=& 2 \exp
\biggl( - \frac{\hbar^2 |a-1|^2 \log p}{ 1+o(1)} \biggr)\\
 &=& 2 p^{-\vafrac{\hbar^2 |a-1|^2}{1+o(1)}}. 
\end{eqnarray*}
\upqed\end{pf}

\begin{lem} \label{lem-norm}
%
\begin{eqnarray}\label{E-trace-norm}
 E \| \hat{\brho} - \brho\|_F^2 &=&  p^{-1/2}
\sum_{j=2}^p E|\hat {\beta}_j
- \beta_j|^2,
\\
\label{E-spectral-norm}
p^{1/2} E \| \hat{\brho} - \brho\|_2 &\leq & \sum
_{j=2}^p E|\hat {\beta}_j -
\beta_j|,
\\
p E \| \hat{\brho} - \brho\|_2^2 &\leq & \sum
_{j=2}^p E\bigl[|\hat {\beta}_j -
\beta_j|^2\bigr] + \Biggl\{ \sum
_{j=2}^p E\bigl[|\hat{\beta}_j -
\beta_j|\bigr] \Biggr\}^2
\nonumber
\\[-8pt]
\label{E-spectral-norm-square}
\\[-8pt]
\nonumber
&& {}- \sum_{j=2}^p \bigl\{E\bigl(|\hat{
\beta}_j - \beta _j|\bigr)\bigr\}^2.
\end{eqnarray}
\end{lem}

\begin{pf}
Since Pauli matrices $\bB_j$ are orthogonal with respect to the
usual Euclidean inner product, with $\|\bB_j\|_F=d^{1/2}$, and $\|\bB
_j\|_2=1$, we have
%
\begin{eqnarray}
\|\hat{\brho} - \brho\|_F^2 &=& \Biggl\| \sum
_{j=2}^p (\hat{\beta}_j -
\beta_j) \bB_j\Biggr \|^2_F\Big/d^2
= \sum_{j=2}^p |\hat{
\beta}_j - \beta _j|^2 \| \bB_j
\|_F^2/d^2
\nonumber
\\[-8pt]
\label{trace-norm}
\\[-8pt]
\nonumber
& =& \sum_{j=2}^p |\hat{
\beta}_j - \beta_j|^2/d,
\\
 p^{1/2} \| \hat{\brho} - \brho\|_2 &=& \Biggl\| \sum
_{j=2}^p (\hat{\beta }_j -
\beta_j) \bB_j \Biggr\|_2 \leq\sum
_{j=2}^p |\hat{\beta}_j - \beta
_j| \| \bB_j \|_2
\nonumber
\\[-8pt]
\label{spectral-norm}
\\[-8pt]
\nonumber
& =& \sum_{j=2}^p |\hat{
\beta}_j - \beta_j|,
\\
p \| \hat{\brho} - \brho\|_2^2  &=& \Biggl\| \sum
_{j=2}^p (\hat{\beta}_j -
\beta_j) \bB_j \Biggr\|_2^2
\nonumber\\
&\leq & \sum_{j=2}^p |\hat{
\beta}_j - \beta_j|^2 \| \bB_j
\|_2^2 + 2 \sum_{i < j }^p
\bigl|(\hat{\beta}_i - \beta_i) (\hat{\beta}_j
- \beta_j)\bigr| \| \bB_i \bB_j
\|_2
\nonumber
\\[-8pt]
\label{spectral-norm-square}
\\[-8pt]
\nonumber
& \leq & \sum_{j=2}^p |\hat{
\beta}_j - \beta_j|^2 \| \bB_j
\|_2^2 + 2 \sum_{i < j }^p
\bigl|(\hat{\beta}_i - \beta_i) (\hat{\beta}_j
- \beta_j)\bigr| \| \bB_i\|_2 \|
\bB_j \|_2
\nonumber
\\
& = &\sum_{j=2}^p |\hat{
\beta}_j - \beta_j|^2 + 2 \sum
_{i < j }^p\bigl|(\hat{\beta}_i -
\beta_i) (\hat{\beta}_j - \beta_j)\bigr|.
\nonumber
\end{eqnarray}

As $N_2, \ldots, N_p$ are independent, $\hat{\beta}_2, \ldots, \hat
{\beta}_p$ are independent. Thus, from
(\ref{trace-norm})--(\ref{spectral-norm-square}) we obtain (\ref
{E-trace-norm})--(\ref{E-spectral-norm}), and
\begin{eqnarray*}
 && p E \| \hat{\brho} - \brho
\|_2^2 \\
&&\qquad\leq  \sum_{j=2}^p
E |\hat {\beta}_j - \beta_j|^2 + 2 \sum
_{i < j }^p E \bigl|(\hat{\beta}_i -
\beta_i) (\hat{\beta}_j - \beta_j)\bigr|
\nonumber
\\
&&\qquad = \sum_{j=2}^p E |\hat{
\beta}_j - \beta_j|^2 + 2 \sum
_{i < j
}^p E |\hat{\beta}_i -
\beta_i| E |\hat{\beta}_j - \beta_j|
\nonumber
\\
&&\qquad= \sum_{j=2}^p E\bigl[|\hat{
\beta}_j - \beta_j|^2\bigr] + \Biggl\{ \sum
_{j=2}^p E\bigl[|\hat{\beta}_j -
\beta_j|\bigr] \Biggr\}^2 - \sum_{j=2}^p
\bigl\{E\bigl(|\hat{\beta}_j - \beta_j|\bigr)\bigr
\}^2.
\end{eqnarray*}
\upqed\end{pf}

\begin{remark}\label{rem5}
Since Pauli matrices $\bB_j$ have eigenvalues $\pm1$, the Schatten
$s$-norm $\| \bB_j\| _{*s}= d^{1/s}$. Similar to
(\ref{spectral-norm})--(\ref{spectral-norm-square}), we obtain that
%
\begin{eqnarray}
\label{s-norm}
 p^{1/2} \| \hat{\brho} - \brho\|_{*s}
&\leq & \sum_{j=2}^p |
\hat{\beta}_j - \beta_j| \| \bB_j
\|_{*s} = d^{1/s} \sum_{j=2}^p
|\hat{\beta}_j - \beta_j|,
\\
 p \| \hat{\brho} - \brho\|_{*s}^2 
&\leq &  d^{2/s} \Biggl[ \sum_{j=2}^p
|\hat{\beta}_j - \beta_j| \Biggr]^2
\nonumber
\\[-8pt]
\label{s-norm-square}
\\[-8pt]
\nonumber
& =&
d^{2/s} \Biggl[ \sum_{j=2}^p |
\hat{\beta}_j - \beta_j|^2 + 2 \sum
_{i < j }^p \bigl|(\hat{\beta}_i -
\beta_i) (\hat{\beta}_j - \beta_j)\bigr| \Biggr].\hspace*{-12pt}
\nonumber
\end{eqnarray}
\end{remark}

\begin{lem} \label{lem-thm3}
%
\begin{eqnarray}
\label{beta-3}  \sum_{j=2}^p E|\hat{
\beta}_j - \beta_j| &\leq &  C_1
\pi_n(d) \varpi^{1-q},
\\
\sum_{j=2}^p \bigl[E|\hat{
\beta}_j - \beta_j|\bigr]^2 &\leq& \sum
_{j=2}^p E\bigl[|\hat{\beta}_j -
\beta_j|^2\bigr] \leq C_2
\pi_n(d) \varpi^{2-q}.
\end{eqnarray}
\end{lem}
\begin{pf}
Using (\ref{threshold1}), we have
\begin{eqnarray*}
&&  E|\hat{\beta}_j - \beta_j|\\
&&\qquad \leq  E \bigl[ \bigl(|
N_j - \beta_j| + \varpi\bigr) 1\bigl(|N_j| \geq
\varpi\bigr) \bigr]+ |\beta_j| P\bigl(|N_j| \leq\varpi \bigr)
\nonumber
\\
&&\qquad \leq  \bigl[E| N_j - \beta_j|^2
P\bigl(|N_j| \geq\varpi\bigr) \bigr]^{1/2} + \varpi
P\bigl(|N_j| \geq\varpi\bigr) + |\beta_j| P\bigl(|N_j|
\leq\varpi\bigr)
\nonumber
\\
&&\qquad \leq   \bigl[ n^{-1} \bigl(1-\beta_j^2\bigr)
P\bigl(|N_j| \geq\varpi\bigr) \bigr]^{1/2} +\varpi
P\bigl(|N_j| \geq\varpi\bigr) + |\beta_j| P\bigl(|N_j|
\leq\varpi \bigr)
\nonumber
\\
&&\qquad\leq  2 \varpi \bigl[P\bigl(|N_j| \geq\varpi\bigr) \bigr]^{1/2} + |
\beta _j| P\bigl(|N_j| \leq\varpi\bigr)
\nonumber
\\
&&\qquad=  2 \varpi \bigl[P\bigl(|N_j| \geq\varpi\bigr) \bigr]^{1/2}
\bigl\{1\bigl(|\beta _j| > a_1 \varpi\bigr)+ 1\bigl(|\beta_j|
\leq a_1 \varpi\bigr)\bigr\}
\nonumber
\\
&&\qquad\quad{}+ |\beta_j| P\bigl(|N_j| \leq\varpi\bigr) \bigl\{1\bigl(|
\beta_j| > a_2 \varpi\bigr)+1\bigl(|\beta_j| \leq
a_2 \varpi\bigr)\bigr\}
\nonumber
\\
&&\qquad\leq   2\varpi1\bigl(|\beta_j| > a_1 \varpi\bigr) +2 \varpi
\bigl[P\bigl(|N_j| \geq\varpi\bigr) \bigr]^{1/2} 1\bigl(|
\beta_j| \leq a_1 \varpi\bigr)
\nonumber
\\
&&\qquad\quad{}+ P\bigl(|N_j|\leq\varpi\bigr) 1\bigl(|\beta_j| > a_2
\varpi\bigr) + |\beta _j| 1\bigl(|\beta_j| \leq a_2
\varpi\bigr),
\end{eqnarray*}
where $a_1$ and $a_2$ are two constants satisfying $a_1 < 1 < a_2$
whose values will be chosen later, and
\begin{eqnarray}
\sum_{j=2}^p E|\hat{
\beta}_j - \beta_j| &\leq &  2\varpi\sum
_{j=2}^p 1\bigl(|\beta_j| >
a_1 \varpi\bigr)\nonumber
\\
&& {}+ 2\varpi\sum_{j=2}^p
\bigl[P\bigl(|N_j| \geq \varpi\bigr) \bigr]^{1/2} 1\bigl(|
\beta_j| \leq a_1 \varpi\bigr)\label{beta-1}
\\
&& {}+ \sum_{j=2}^p P\bigl(|N_j|
\leq\varpi\bigr) 1\bigl(|\beta _j| > a_2 \varpi\bigr) + \sum
_{j=2}^p |\beta_j| 1\bigl(|
\beta_j| \leq a_2 \varpi\bigr).
\nonumber
\end{eqnarray}
Similarly,
\begin{eqnarray*}
 && \bigl[E\bigl(|\hat{\beta}_j - \beta_j|\bigr)
\bigr]^2\\
 &&\qquad \leq   E\bigl[|\hat{\beta}_j - \beta
_j|^2\bigr]
\nonumber
\\
& &\qquad\leq   E\bigl[ 2\bigl(| N_j - \beta_j|^2
+ \varpi^2 \bigr) 1\bigl(|N_j| \geq\varpi\bigr)\bigr] + |
\beta_j|^2 P\bigl(|N_j| \leq\varpi\bigr)
\nonumber
\\
&&\qquad \leq  2 \bigl[E| N_j - \beta_j|^4
P\bigl(|N_j| \geq\varpi\bigr) \bigr]^{1/2}\\
&&\qquad\quad {}+2 \varpi ^2
P\bigl(|N_j| \geq\varpi\bigr) + |\beta_j|^2
P\bigl(|N_j| \leq\varpi\bigr)
\nonumber
\\
&&\qquad\leq   c \varpi^2 \bigl[P\bigl(|N_j| \geq\varpi\bigr)
\bigr]^{1/2} + |\beta_j|^2 P\bigl(|N_j|
\leq\varpi\bigr)
\nonumber
\\
&&\qquad =  c \varpi^2 \bigl[P\bigl(|N_j| \geq\varpi\bigr)
\bigr]^{1/2} \bigl\{1\bigl(|\beta_j| > a_1 \varpi\bigr)+
1\bigl(|\beta_j| \leq a_1 \varpi\bigr)\bigr\}
\nonumber
\\
&& \qquad\quad{}+ |\beta_j|^2 P\bigl(|N_j| \leq\varpi\bigr)
\bigl[1\bigl(|\beta_j| > a_2 \varpi\bigr) + 1\bigl(|
\beta_j| \leq a_2 \varpi\bigr)\bigr]
\nonumber
\\
&&\qquad \leq   c \varpi^2 1\bigl(|\beta_j| > a_1
\varpi\bigr) + c \varpi^2 \bigl[P\bigl(|N_j| \geq\varpi\bigr)
\bigr]^{1/2} 1\bigl(|\beta_j| \leq a_1 \varpi\bigr)
\nonumber
\\
&& \qquad\quad{}+ P\bigl(|N_j|\leq\varpi\bigr) 1\bigl(|\beta_j| > a_2
\varpi\bigr) + |\beta_j|^2 1\bigl(|\beta_j| \leq
a_2 \varpi\bigr),
\end{eqnarray*}
and
%
\begin{eqnarray}
&&\sum_{j=2}^p E\bigl[|
\hat{\beta}_j - \beta_j|^2\bigr]\nonumber\\
 &&\qquad \leq   c
\varpi^2 \sum_{j=2}^p 1\bigl(|
\beta_j| > a_1 \varpi\bigr)\nonumber
\\[-8pt]
\label{beta-2}
\\[-8pt]
\nonumber
&&\qquad\quad {}+c \varpi^2 \sum_{j=2}^p
\bigl[P\bigl(|N_j| \geq\varpi \bigr) \bigr]^{1/2} 1\bigl(|
\beta_j| \leq a_1 \varpi\bigr)
\\
&& \qquad\quad{}+ \sum_{j=2}^p P\bigl(|N_j|
\leq\varpi\bigr) 1\bigl(|\beta_j| > a_2 \varpi\bigr) + \sum
_{j=2}^p |\beta_j|^2 1\bigl(|
\beta_j| \leq a_2 \varpi \bigr).
\nonumber
\end{eqnarray}
By Lemma~\ref{lem-sparse}, we have\vspace*{-2pt}
%
\begin{eqnarray}
\label{beta-sparse-1} &&\sum_{j=2}^{p} |\beta_j|
1\bigl(|\beta_j| \leq a_2 \varpi\bigr)  \leq   a_2^{1-q}
\pi_n(d) \varpi^{1-q},
\\
&& \sum_{j=2}^{p} |\beta_j|^2
1\bigl(|\beta_j| \leq a_2 \varpi\bigr)
\nonumber
\\[-8pt]
\label{beta-sparse-11}
\\[-8pt]
\nonumber
&&\qquad \leq  (a _2 \varpi)^{2-q} \sum
_{j=2}^{p} |\beta_j|^q 1\bigl(|
\beta_j| \leq a_2 \varpi\bigr) \leq a_2^{2-q}
\pi_n(d) \varpi^{2-q},
\\
\label{beta-sparse-2}
&& \varpi \sum_{j=2}^p 1\bigl(|
\beta_j| \geq a_1 \varpi\bigr) \leq a_1^{-q}\pi_n(d)
\varpi^{1-q}.
\end{eqnarray}
On the other hand,\vspace*{-2pt}
%
\begin{eqnarray}
&& \sum_{j=2}^p
P\bigl(|N_j|\leq\varpi\bigr) 1\bigl(|\beta_j| > a_2
\varpi\bigr)\nonumber
\\
&&\qquad \leq\sum_j P( -\varpi- \beta_j
\leq N_j - \beta_j \leq\varpi- \beta_j) 1\bigl(|
\beta_j| > a_2 \varpi\bigr)
\nonumber
\\
\label{beta-sparse-3}
&& \qquad\leq\sum_{j=2}^p \bigl[P\bigl(
N_j - \beta_j \leq-|a_2-1| \varpi\bigr) +
P\bigl(N_j - \beta_j \geq|a_2-1| \varpi\bigr)\bigr]
\\
&& \qquad\leq4 p^{1- \vafrac{\hbar^2 |a_2-1|^2}{1+o(1)}} = 4 p^{-1 -
(2-q)/(2 c_0)} \leq4 p^{-1}
n^{-(q-2)/2}\nonumber\\
&&\qquad = o\bigl(\pi_n(d) \varpi^{2-q}\bigr),
\nonumber
\end{eqnarray}
where the third inequality is from Lemma~\ref{lem-tail}, the first
equality is due the fact that we take
$a_2 = 1 + \{2 + (2 - q)/(2 c_0)\}^{1/2} (1+o(1))^{1/2}/\hbar$ so that
$\hbar^2 (1-a_2)^2/ (1+o(1))=2 + (2-q)/(2 c_0)$,
and $c_0$ is the constant in assumption $p \geq n^{c_0}$.
Finally, we can show\vspace*{-2pt}
%
\begin{eqnarray}
\nonumber
&& \varpi \sum_{j=2}^p
\bigl[ P\bigl(|N_j| \geq\varpi\bigr) \bigr]^{1/2} 1\bigl(|
\beta_j| \leq a_1 \varpi\bigr)
\\
&& \qquad\leq\varpi \sum_{j=2}^p \bigl[P(
N_j - \beta_j \leq- \varpi- \beta _j)
\nonumber
\\
\label{beta-sparse-4}
&& \quad\qquad{}+ P(N_j - \beta_j \geq\varpi- \beta_j
)\bigr]^{1/2} 1\bigl(|\beta_j| \leq a_1 \varpi\bigr)
\\
&& \qquad\leq\varpi \sum_{j=2}^p \bigl[P\bigl(
N_j - \beta_j \leq- |1-a_1| \varpi \bigr) +
P\bigl(N_j - \beta_j \geq|1-a_1| \varpi\bigr)
\bigr]^{1/2}
\nonumber
\\
&&\qquad \leq2\varpi p^{1- \hbar^2 (1-a_1)^2/(2(1+o(1)))} =2\varpi p^{-1} = o\bigl(
\pi_n(d) \varpi^{1-q}\bigr),
\nonumber
\end{eqnarray}
where the third inequality is from Lemma~\ref{lem-tail}, and the first
equality is due to the fact that we take
$a_1 = 1 - 2 (1+o(1))^{1/2}/\hbar$ so that $\hbar^2 (1-a_1)^2=4$.
Plugging (\ref{beta-sparse-1})--(\ref{beta-sparse-4}) into (\ref{beta-1}) and (\ref{beta-2}),
we prove the lemma.
\end{pf}

\begin{pf*}{Proof of Theorem~\protect\ref{thm-3}}
Combining Lemma~\ref
{lem-thm3} and
(\ref{E-trace-norm})--(\ref{E-spectral-norm}) in Lemma~\ref
{lem-norm}, we easily obtain
\begin{eqnarray*}
E\bigl[\| \hat{\brho} - \brho\|_2\bigr] &\leq&  C_1
\frac{\pi_n(d)}{p^{1/2}} \biggl( \frac{ \log p }{n} \biggr)^{\vfrac{1-q}{2}},
\\
E\bigl[\| \hat{\brho} - \brho\|_F^2\bigr]  &\leq &
C_0 \pi_n(d) \frac{1}{d} \biggl(
\frac{ \log p }{n} \biggr)^{1-q/2}.
\end{eqnarray*}
%
Using Lemma~\ref{lem-thm3} and (\ref{E-spectral-norm-square}) in
Lemma~\ref{lem-norm}, we conclude
%
\begin{eqnarray}
 E\bigl[\| \hat{\brho} - \brho\|_2^2
\bigr]  &\leq &  C_2 \biggl[ \pi^2_n(d)
\frac{1}{p} \biggl( \frac{ \log p }{n} \biggr)^{1-q} +
\pi_n(d) \frac{1}{p} \biggl( \frac{ \log p }{n}
\biggr)^{1-q/2} \biggr]
\nonumber
\\[-8pt]
\label{equation-spectral-norm-square}
\\[-8pt]
\nonumber
& \leq &  C \frac{\pi^2_n(d)}{d^2} \biggl( \frac{ \log p }{n} \biggr)^{1-q},
\end{eqnarray}
where the last inequality is due to the fact that the first term on the
right-hand side of (\ref{equation-spectral-norm-square}) dominates
its second term.
\end{pf*}

\begin{pf*}{Proof of Proposition~\protect\ref{schatten.prop}}
Applying the
Lyapunov's moment inequality to the Schatten $s$-norm, we have for $s
\in[1, 2]$ and $\brho\in\Theta$,
\begin{eqnarray*}
E \bigl[ \| \hat{\brho}- \brho\|_{*s} ^2 \bigr] &\leq&
d^{-1 +
2/s} E \bigl[ \| \hat{\brho} -\brho\|_{*2} ^2
\bigr]
\\
&=& d^{-1 + 2/s} E \bigl[ \| \hat{\brho} -\brho\|_{F}
^2 \bigr]
\\
&\leq& c_1 \pi_n (d)
d^{-2 + 2/s} \biggl(\frac{\log d}{n} \biggr) ^{1-q/2},
\end{eqnarray*}
where the last inequality is due to Theorem~\ref{thm-3}.
On the other hand, applying H\"older's inequality by interpolating
between Schatten $s$-norms with $s= 2$ and $s=\infty$,
we obtain for $s \in[2, \infty]$ and $\brho\in\Theta$, 
\begin{eqnarray*}
E \bigl[ \| \hat{\brho}- \brho\|_{*s}^2 \bigr] &\leq& E
\bigl[ \| \hat{\brho}- \brho\|_{*2}^{4/s} \| \hat{\brho}-
\brho\| _{*\infty} ^{2-4/s} \bigr]
\\
&\leq& \bigl[ E \| \hat{
\brho}- \brho\|_{*2}^{2} \bigr] ^{2/s} \bigl[E \|
\hat{\brho}- \brho\|_{*\infty} ^{2} \bigr] ^{1-2/s}
\\
&
\leq& c_7 \pi_n ^{2-2/s}(d) d^{-2+2/s}
\biggl(\frac{\log d}{n} \biggr) ^{1-q+q/s},
\end{eqnarray*}
where the last inequality is due to Theorem~\ref{thm-3}, and $c_7 =
c_1^{(s-2)/s} c_2^{2/s}$. The result follows by combining the above two
inequalities together.
\end{pf*}

\subsection{Proofs of Theorems \texorpdfstring{\protect\ref{thm-4}}{2} and \texorpdfstring{\protect\ref{thm-5}}{3}: Lower bound}
\label{lower}
\mbox{}
\begin{pf*}{Proof of Theorem~\protect\ref{thm-4} \textup{(for the lower bound under the
spectral norm)}}
We first define a subset of the parameter space $\Theta$. It will be shown
later that the risk upper bound under the spectral norm is sharp up to a
constant factor, when the parameter space is sufficiently sparse.
Consider a
subset of the Pauli basis, $ \{ \bolds{\sigma}_{l_{1}}\otimes
\bolds{\sigma}_{l_{2}}\otimes\cdots\otimes\bolds{\sigma}_{l_{b}}
\} $,
where $\bolds{\sigma}_{l_{1}}=\bolds{\sigma}_{0}$ or $\bolds
{\sigma}
_{3}$. Its cardinality is $d=2^{b}=p^{1/2}$. Denote each element of the
subset by $\mathbf{B}_{j}$, $j=1,2,\ldots,d$, and let $\mathbf{B}_{1}=
\mathbf{I}_{d}$. We will define each element of $\Theta$ as a linear
combination of $\mathbf{B}_{j}$. Let $\gamma_{j}\in \{
0,1 \} $, $
j\in \{ 1,2,\ldots,d \} $, and denote $\eta=\sum_{j}\gamma
_{j}=\llVert \bolds{\gamma}\rrVert _{0}$. The value of
$\eta$ is
either $0$ or~$K$, where $K$ is the largest integer less than or equal
to $
\pi_{n} ( d  ) / ( \frac{\log p}{n} ) ^{q/2}$.
By assumption (\ref{cond-pi}),
we have
%
\begin{equation}\label{K}
1\leq K=O \bigl( d^{v} \bigr)\qquad \mbox{with }v<1/2.
\end{equation}
Let $\varepsilon^{2}= ( 1-2v ) /4$ and set $a=\varepsilon
\sqrt{
\frac{\log p}{n}}$. Now we are ready to define $\Theta$,
%
\begin{equation}\label{Subspace}
\Theta= \Biggl\{ \bolds{\rho} ( \bolds{\gamma} ) :\bolds{\rho} ( \bolds{\gamma}
) =\frac{\mathbf{I}_{d}}{d}+a\sum_{j=2}^{d}
\gamma_{j}\frac{\mathbf{B}_{j}}{d},\mbox{ and }\eta=0\mbox{ or }K \Biggr\}.
\end{equation}
Note that $\Theta$ is a subset of the parameter space, since
\[
\sum_{j=2}^{d} ( a\gamma_{j} )
^{q}\leq Ka^{q}\leq \varepsilon ^{q}
\pi_{n} ( d ) \leq\pi_{n} ( d ),
\]
and its cardinality is $1+ {d-1 \choose K}$.

We need to show that
\[
\inf_{\bolds{\hat{\brho}}}\sup_{\Theta}E\llVert\hat{ \bolds{\brho}}-\bolds{\brho} \rrVert _{2}^{2}\gtrsim\pi_{n}^{2}
( d ) \frac
{1}{p} \biggl( \frac{\log p}{n} \biggr) ^{1-q}.
\]
Note that for each element in $\Theta$, its first entry $\bolds{\rho
}_{11}$ may take the form
$1/d+a\sum_{j=2}^{d}\gamma_{j}/d=1/d+ (a/d) \eta$. It can be shown that
\[
\inf_{\hat{\bolds{\brho}}}\sup_{\Theta}E\llVert \hat{\bolds{\brho}}-\bolds{\brho} \rrVert _{2}^{2}\geq\inf_{\hat{\rho}_{11}}
\sup_{\Theta
}E ( \hat{ \rho}_{11}\mathbf{-}
\rho_{11} ) ^{2}\geq\frac
{a^{2}}{d^{2}}\inf
_{\hat{
\eta}}\sup_{\Theta}E ( \hat{\eta}\mathbf{-}
\eta ) ^{2}.
\]
It is then enough to show that
%
\begin{equation}\label{Klowerbd}
\inf_{\hat{\eta}}\sup_{\Theta}E ( \hat{\eta}
\mathbf{-} \eta ) ^{2}\gtrsim K^{2},
\end{equation}
which immediately implies
\[
\inf_{\bolds{\hat{\brho}}}\sup_{\Theta}E\llVert \hat{\bolds{\brho}}-\bolds{\brho} \rrVert _{2}^{2}\gtrsim K^{2}
\frac{a^{2}}{d^{2}}\gtrsim\pi _{n}^{2} ( d )
\frac{1}{p} \biggl( \frac{\log p}{n} \biggr) ^{1-q}.
\]

We prove equation (\ref{Klowerbd}) by applying Le Cam's lemma. From
observations $N_{j}$, $j=2,\ldots,d$, we define $\tilde{N}_{j}=n (
N_{j}+1 ) /2$,
which is $\operatorname{Binomial} ( n,\frac{1+a\gamma_{j}}{2} ) $. Let
$\mathbb{P}
_{\bolds{\gamma}}$ be the joint distribution of independent random
variables $\tilde{N}_{2},\tilde{N}_{3},\ldots,\tilde{N}_{d}$. The
cardinality of $ \{ \mathbb{P}_{\bolds{\gamma}} \} $ is
$1+ {d-1 \choose K}$.
For two probability measures $\mathbb{P}$ and $\mathbb{Q}$ with
density $f$ and $g$ with respect to any common dominating measure $\mu$,
write the total variation affinity $\Vert\mathbb{P}\wedge\mathbb
{Q}\Vert
=\int f\wedge g \,d\mu$, and the Chi-square distance $\chi^{2} (
\mathbb{P
},\mathbb{Q} ) =\int\frac{g^{2}}{f}-1$. Define
\[
\bar{ \mathbb{P}} = \pmatrix{d-1 \cr K}^{-1}\sum
_{\llVert \bolds
{\gamma}
\rrVert _{0}=K}\mathbb{P}_{\bolds{\gamma}}.
\]
The following lemma is a direct consequence of Le Cam's lemma
[cf. \citet{LeC73} and \citet{Yu97}].
\end{pf*}

\begin{lem}
\label{LeCam}Let $\hat{\eta}$ be any estimator of $\eta$ based on
an observation
from a distribution in the collection $ \{ \mathbb{P}_{\bolds
{\gamma}
} \} $, then
\[
\inf_{\hat{k}}\sup_{\Theta}E ( \hat{\eta}
\mathbf{-}\eta ) ^{2}\geq\frac{
1}{4}\llVert
\mathbb{P}_{\mathbf{0}}\wedge\bar{\mathbb {P}}\rrVert ^{2}\cdot
K^{2}.
\]
\end{lem}

We will show that there is a constant $c>0$ such that
%
\begin{equation}\label{aff}
\llVert \mathbb{P}_{\mathbf{0}}\wedge\bar{\mathbb{P}}\rrVert \geq C,
\end{equation}
which, together with Lemma~\ref{LeCam}, immediately imply equation
(\ref{Klowerbd}).

\begin{lem}
Under conditions (\ref{K}) and (\ref{Subspace}), we have
\[
\inf_{\bolds{\hat{\brho}}}\sup_{\Theta}E ( \hat{\eta }
\mathbf{-}\eta ) ^{2}\gtrsim K^{2},
\]
which implies
\[
\inf_{\bolds{\hat{\brho}}}\sup_{\Theta}E\llVert \hat{\bolds{\brho}}-\bolds{\brho} \rrVert _{2}^{2}\gtrsim\pi_{n}^{2}
( d ) \frac
{1}{p} \biggl( \frac{\log p}{n} \biggr) ^{1-q}.
\]
\end{lem}

\begin{pf}
It is enough to show that
\[
\chi^{2} ( \mathbb{P}_{\mathbf{0}},\bar{\mathbb{P}} )
\rightarrow0,
\]
which implies $\llVert \mathbb{P}_{\mathbf{0}}-\bar{\mathbb{P}}
\rrVert _{\mathrm{TV}}\rightarrow0$, then we have $\llVert \mathbb{P}_{
\mathbf{0}}\wedge\bar{\mathbb{P}}\rrVert \rightarrow1$. Let
$J (
\bolds{\gamma},\bolds{\gamma}^{{\prime}} ) $ denote the number of
overlapping nonzero coordinates between $\bolds{\gamma}$ and
$\bolds{
\gamma}^{{\prime}}$. Note that
\begin{eqnarray*}
\chi^{2} ( \mathbb{P}_{\mathbf{0}},\bar{\mathbb{P}} ) &=&\int
\frac{ ( d \bar{ \mathbb{P} } ) ^{2}}{d\mathbb
{P}_{\mathbf{0}}}-1
\\
&=& \pmatrix{d-1 \cr K}^{-2}\sum_{0\leq j\leq K}\sum
_{J ( \bolds{\gamma},\bolds{\gamma}^{{\prime}} ) =j} \biggl( \int\frac{d\mathbb
{P}_{\bolds{
\gamma}}\cdot d\mathbb{P}_{\bolds{\gamma}^{{\prime}}}}{d\mathbb{P}_{
\mathbf{0}}}-1 \biggr).
\end{eqnarray*}

When $J ( \bolds{\gamma},\bolds{\gamma}^{{\prime}} ) =j$, we have
\begin{eqnarray*}
\int\frac{d\mathbb{P}_{\bolds{\gamma}}\cdot d\mathbb
{P}_{\bolds{\gamma
}^{{\prime}}}}{d\mathbb{P}_{\mathbf{0}}}
&=& \Biggl( \sum_{l=0}^{n} \biggl[ \pmatrix{n
\cr l} \frac{1}{2^{l}}\frac
{1}{2^{n-l}} \cdot ( 1+a ) ^{2l} (
1-a ) ^{2n-2l} \biggr] \Biggr) ^{j}
\\
&= &\Biggl( \sum_{l=0}^{n} \biggl[ \pmatrix{n \cr l} \biggl( \frac{ (
1+a )
^{2}}{2} \biggr) ^{l} \biggl(
\frac{ ( 1-a ) ^{2}}{2} \biggr) ^{n-l} \biggr] \Biggr) ^{j}
\\
&=& \biggl( \frac{ ( 1+a ) ^{2}}{2}+\frac{ ( 1-a
) ^{2}}{2} \biggr) ^{nj}
\\
&=& \bigl( 1+a^{2} \bigr) ^{nj}\leq\exp \bigl(
na^{2}j \bigr),
\end{eqnarray*}
which implies
\begin{eqnarray*}
\chi^{2} ( \mathbb{P}_{\mathbf{0}},\bar{\mathbb{P}} ) &\leq& \pmatrix{d-1
\cr K}^{-2}\sum_{0\leq j\leq K}\sum
_{J ( \bolds{\gamma},\bolds{\gamma}
^{{\prime}} ) =j} \bigl( \exp \bigl( na^{2}j \bigr) -1 \bigr)
\\
&\leq& \pmatrix{d-1 \cr K}^{-2}\sum_{1\leq j\leq K}\sum
_{J (
\bolds{\gamma},\bolds{\gamma}^{{\prime}} ) =j}\exp \bigl( na^{2}j \bigr)
\\
&=&\sum_{1\leq j\leq K}\frac{ {K \choose j} {d-1-K \choose K-j}
}{ {d-1 \choose K} } d^{2\varepsilon^{2}j}.
\end{eqnarray*}

Since
\begin{eqnarray*}
\frac{ {K \choose j} {d-1-K \choose K-j} }{ {d-1 \choose K} } 
&=&\frac{ [ K\cdot
\ldots\cdot ( K-j+1 )  ] ^{2}\cdot (
d-1-K ) \cdot
\ldots\cdot ( d-2K+j ) }{j!\cdot ( d-1 ) \cdot
\ldots
\cdot ( d-K ) }
\\
&\leq&\frac{K^{2j} ( d-1-K ) ^{K-j}}{ ( d-K )
^{K}}\leq \biggl( \frac{K^{2}}{d-K} \biggr) ^{j},
\end{eqnarray*}
and $\varepsilon^{2}= ( 1-2v ) /4$, we then have
\begin{eqnarray*}
&& \chi^{2} ( \mathbb{P}_{\mathbf{0}},\bar{\mathbb{P}} ) \leq \sum
_{1\leq j\leq K} \biggl[ \frac{K^{2}}{d-K}d^{2\varepsilon
^{2}}
\biggr] ^{j}
\leq\sum_{1\leq j\leq K} \biggl[ \frac{d^{2v+ ( 1-2v ) /2}}{d-K} \biggr]
^{j}\rightarrow0.
\end{eqnarray*}
\upqed\end{pf}

\begin{pf*}{Proof of Theorem~\protect\ref{thm-5} \textup{(for the lower bound under the
Frobenius norm)}}
Recall that $\Theta$ is the collection of density matrices such that
\[
\brho= \frac{1}{d} \Biggl( \bI_d + \sum
_{j=2}^p \beta_j \bB_j
\Biggr),
\]
where
\[
\sum_{j=2}^p |\beta_j|^q
\leq\pi_n(p).
\]

Apply Assouad's lemma, and we show below that
\[
\inf_{\check{\brho}} \sup_{\brho\in\Theta} E\bigl[\| \check{
\brho} - \brho\|_F^2\bigr] \geq C \pi_n(p)
\frac{1}{d} \biggl( \frac{ \log p }{n} \biggr)^{1-q/2},
\]
where $\check{\brho}$ denotes any estimator of $\brho$ based on
measurement data $N_2, \ldots, N_p$, and
$C$ is a constant free of $n$ and $p$.

To this end, it suffices to construct a collection of $M+1$ density
matrices $\{\brho_0=\bI_d/d, \brho_1, \ldots, \brho_M\} \subset
\Theta$ such that
(i) for any distinct $k$ and $k_0$,
\[
\| \brho_k - \brho_{k_0}\|_F^2
\geq C_1 \pi_n(p) \frac{1}{d} \biggl(
\frac{ \log p }{n} \biggr)^{1-q/2},
\]
where $C_1$ is a constant;
(ii) there exists a constant $0<C_2<1/8$ such that
\[
\frac{1}{M} \sum_{k=1}^M
D_{\mathrm{KL}}(P_{\brho_k}, P_{\brho_0}) \leq C_2 \log
M,
\]
where $D_{\mathrm{KL}}$ denotes the Kullback--Leibler divergence.

By the Gilbert--Varshamov bound [cf. \citet{NieChu00}], we
have that for any $h<p/8$, there exist $M$ binary vectors
$\bgamma_k = (\gamma_{k2}, \ldots, \gamma_{kp})^\prime\in\{0,1\}
^{p-1}$, $k=1, \ldots, M$,
such that (i) $\|\bgamma_k\|_{1} = \sum_{j=2}^p |\gamma_{kj}| = h$,
(ii) $\| \bgamma_k - \bgamma_{k_0} \|_{1} =
\sum_{j=2}^p |\gamma_{kj}-\gamma_{k_0 j}| \geq h/2$, and (iii) $\log
M > 0.233  h  \log(p/h)$.
Let
\[
\brho_k = \frac{1}{d} \Biggl( \bI_d + \epsilon
\sum_{j=2}^p \gamma _{kj}
\bB_j \Biggr),
\]
where
\[
\epsilon= C_3 \biggl( \frac{\pi_n(p)}{h} \biggr)^{1/q}.
\]
Since $\sum_{j=2}^p |\epsilon\gamma_{kj}|^q = \epsilon^q h = C_3
\pi_n(p)$, $\brho_k \in\Theta$ whenever $C_3 \leq1$. Moreover,
\[
d \| \brho_k - \brho_{k_0}\|_F^2
= \epsilon^2 \| \bgamma_k - \bgamma_{k_0}
\|_{1} \geq \frac{\epsilon^2  h}{4}.
\]
On the other hand,
\begin{eqnarray*}
D_{\mathrm{KL}}(P_{\brho_k}, P_{\brho_0}) &=& h D_{\mathrm{KL}}
\biggl(\operatorname{Bin} \biggl(n, \frac{1+\epsilon}{2} \biggr), \operatorname{Bin} \biggl(n, \frac{1}{2}
\biggr) \biggr)
\cr
&=& h n \frac{\epsilon}{2} \log\frac{1/2+\epsilon}{1/2-\epsilon} \leq
C_4 h n \epsilon^2.
\end{eqnarray*}
%
Now the lower bound can be established by taking
\[
h = \pi_n(p) \biggl( \frac{\log p}{n} \biggr)^{-q/2},
\]
and then
\begin{eqnarray*}
\epsilon&=&  C_3 \biggl( \frac{\log p}{n} \biggr)^{1/2},\qquad
\frac
{\epsilon^2  h}{4} = C_3 \pi_n(p) \biggl(
\frac{\log p}{n} \biggr)^{1-q/2},
\\
C_4 h n \epsilon^2  &=&  C_4 h \log p,\qquad h
\log(p/h) = h \log p - h \log h,
\\
\log h  &\sim & \log\pi_n(p) +
\frac{q}{2} \log n - \frac{q}{2} \log \log p,
\end{eqnarray*}
which are allowed by the assumption $\log\pi_n(p) + \frac{q}{2} \log
n < v^\prime\log p$ for $v^\prime<1$.
\end{pf*}






\printaddresses
\end{document}